\documentclass[12pt]{article}%
\usepackage{amssymb}
\usepackage{amsmath}
\usepackage{graphics}
\usepackage{epsfig}
\usepackage[latin1]{inputenc}%
\usepackage{amsfonts}%
\usepackage{graphicx}
\providecommand{\U}[1]{\protect\rule{.1in}{.1in}}
\DeclareMathOperator{\argsup}{\argsup}\DeclareMathOperator{\arginf}{\arginf} 
\newtheorem{definition}{Definition}
\newtheorem{theorem}{Theorem}

\numberwithin{equation}{section}

\begin{document}

\begin{center}
{\Large \textsc{Decomposable Pseudodistances and Applications\newline%
\vspace{2mm} in Statistical Estimation }}
\end{center}

\vspace{4mm}

\begin{center}
\textsc{Michel Broniatowski\footnote{Laboratoire de Statistique
Théorique et Appliquée, Université Paris 6, Paris, France, e-mail:
michel.broniatowski@upmc.fr }\ \ \ \ \ Aida
Toma\footnote{Mathematics Department, Academy of Economic Studies
and \textit{Gh. Mihoc - C. Iacob} Institute of Mathematical
Statistics and Applied
Mathematics, Bucharest, Romania, e-mail: aida$_{-}$toma@yahoo.com}%
\ \ \ \ \ Igor Vajda}$^{\dag}$\textsc{\footnote{Institute of Information
Theory and Automation, Academy of Sciences of the Czech Republic, Prague,
Czech Republic.} }
\end{center}
\vspace{4mm}

\begin{abstract}
The aim of this paper is to introduce new statistical criterions
for estimation, suitable for inference in models with common
continuous support. This proposal is in the direct line of a
renewed interest for divergence based inference tools imbedding
the most classical ones, such as maximum likelihood, Chi-square or
Kullback Leibler. General pseudodistances with decomposable
structure are considered, they allowing to define minimum
pseudodistance estimators, without using nonparametric density
estimators. A special class of pseudodistances indexed by
$\alpha>0$, leading for $\alpha\downarrow0$ to the Kulback Leibler
divergence, is presented in detail. Corresponding estimation
criteria are developed and asymptotic properties are studied. The
estimation method is then extended to regression models. Finally,
some examples based on Monte Carlo simulations are discussed.

\end{abstract}

\section{Introduction}

In parametric estimation, minimum divergence methods, i.e. methods which
estimate the parameter by minimizing an estimate of some divergence between
the assumed model density and the true density underlying the data, have been
extensively studied (see Pardo (2005) and references herein). Generally, in
continuous models, the minimum divergence methods have the drawback that it is
necessary to use some nonparametric density estimator. In order to remove this
drawback, some proposals have been made in literature. Among them, we recall
the minimum density power divergence method introduced by Basu et al. (1998),
and a minimum divergence method based on duality arguments, independently
proposed by Liese and Vajda (2006) and Broniatowsi and Keziou (2009). The
results obtained in the present paper follow this line of research.

Without referring to all properties of the divergence criterions, we mainly
quote their information processing property, i.e. the complete invariance with
respect to the statistically sufficient transformations of the observation
space. This property is useful but probably not unavoidable in the minimum
divergence estimation based on similarity between theoretical and empirical
distributions. In this paper we admit general pseudodistances which may not
satisfy the information processing property. The definition of the
pseudodistance, which is at the start of this work, pertains to the
willingness to define a simple frame including all commonly used statistical
criterions, from maximum likelihood to the L$_{2}$ norm. Such a description is
provided in Broniatowski and Vajda (2009). In the present paper we define a
class of pseudodistances indexed by $\alpha>0$, leading for $\alpha
\downarrow0$ to the Kulback Leibler divergence. The peculiar features of these
pseudodistances recommend it as an appealing competing choice for defining
estimation criteria. We argue that by defining and studying minimum
pseudodistances estimators for classical parametric models, respectively for
regression models. We present such tools for inference with a special
attention to limit properties and robustness, in a similar spirit as in Toma
and Broniatowski (2011).

The outline of the paper is as follows. Section 2 introduces decomposable
pseudodistances and define minimum pseudodistances estimators. Section 3
presents a special class of minimum pseudodistances estimators. For these
estimators we study invariance properties, consistency, asymptotic normality
and robustness. The estimation method is applied to linear models for which
asymptotic and robustness properties are derived. These results are presented
in Section 4. Finally, in order to illustrate the performance of the proposed
method in finite samples, we give some examples based on Monte Carlo simulations.

\section{Decomposable pseudodistances and estimators}

We will consider inference in continuous parametric families, since this is
the interesting and complex case with respect to the case of models with
finite or countable support. Hence $\mathcal{P}$ is a parametric model with
euclidian parameter space $\Theta$ and we assume that all the probability
measures $P_{\theta}$ in $\mathcal{P}$ share the same support, which is
included in $\mathbb{R}^{d}$. Every $P_{\theta}$ has a density $p_{\theta}$
with respect to the Lebesgue measure.

We denote by ${\mathcal{P}}_{\text{emp}}$ the class of probability measures
induced by samples, namely the class of all probability measures
\[
P_{n}:=\frac{1}{n}\sum_{i=1}^{n}\delta_{X_{i}},
\]
where $X_{1},\dots,X_{n}$ is sampled according to a distribution on $\left(
\mathbb{R}^{d},\mathcal{B}\left(  \mathbb{R}^{d}\right)  \right) $, not
necessarily in $\mathcal{P}.$ In addition to the previous notation it is
useful to introduce a family of measures $\mathcal{P}_{0}$ associated to
distributions generating the data when studying robustness properties. Often,
such a measure is a mixture of some element in $\mathcal{P}$ with a Dirac
measure at some point $x$ in $\mathbb{R}^{d}$. We also define $\mathcal{P}%
^{+}:=\mathcal{P\cup}\mathcal{P}_{0}$.

\begin{definition}
We say that $\mathfrak{D}:\mathcal{P}\otimes\mathcal{P}^{+}\mapsto\mathbb{R}$
is a pseudodistance between a probability measures $P\in\mathcal{P}%
=\{P_{\theta}:\theta\in\Theta\}$ and $Q\in\mathcal{P}^{+}$ if $\mathfrak{D}%
(P_{\theta},Q)\geq0$, for all $\theta\in\Theta$ and $Q\in\mathcal{P}^{+}$ and
$\mathfrak{D}(P_{\theta},P_{{\tilde{\theta}}})=0$ if and only if
$\theta={\tilde{\theta}}$.
\end{definition}

\begin{definition}
A pseudodistance $\mathfrak{D}$ on $\mathcal{P}\otimes\mathcal{P}^{+}$ is
called $decomposable$ if there exist functionals $\mathfrak{D}^{0}%
:\mathcal{P}\mapsto\mathbb{R}$, $\mathfrak{D}^{1}:\mathcal{P}^{+}%
\mapsto\mathbb{R}$ and measurable mappings
\begin{equation}
\rho_{\theta}:\mathbb{R}^{d}\mapsto\mathbb{R},\text{ \ \ \ }\theta\in
\Theta\label{24aa}%
\end{equation}
\ such that for all $\theta\in\Theta$\ and $Q\in\mathcal{P}^{+}$ the
expectations $\int\rho_{\theta}\mathrm{d}Q$\ exist and
\begin{equation}
\mathfrak{D}(P_{\theta},Q)=\mathfrak{D}^{0}(P_{\theta})+\mathfrak{D}%
^{1}(Q)+\int\rho_{\theta}\mathrm{d}Q.\label{26a}%
\end{equation}

\end{definition}

A known class of pseudodistances is that introduced by Basu et al. (1998) and
called the class of power divergences. This class corresponds to
\begin{equation}
\label{Basu}\mathfrak{D}(P_{\theta},Q)=\int\left\{ p_{\theta}^{\alpha
+1}-\left( 1+\frac{1}{\alpha}\right) p_{\theta}^{\alpha}q+\frac{1}{\alpha
}q^{\alpha+1}\right\} \mathrm{d}\lambda\;\;\text{for}\;\;\alpha>0.
\end{equation}
Note that the pseudodistances (\ref{Basu}) are decomposable in the sense
(\ref{26a}) with
\begin{equation}
\mathfrak{D}^{0}(P_{\theta})=\int p_{\theta}^{\alpha+1}\mathrm{d}%
\lambda,\;\mathfrak{D}^{1}(Q)=\frac{1}{\alpha}\int q^{\alpha+1}\mathrm{d}%
\lambda\;\text{and}\;\rho_{\theta}=-\left( 1+\frac{1}{\alpha}\right)
p_{\theta}^{\alpha}.
\end{equation}

In the next section, we introduce a new class of pseudodistances
from which a new statistical criterion for inference is deduced.

\begin{definition}
We say that a functional $T_{\mathfrak{D}}:\mathcal{Q}\mapsto\Theta$ for
$\mathcal{Q}=\mathcal{P}^{+}\cup\mathcal{P}_{\text{emp}}$ defines a minimum
pseudodistance estimator (briefly, $\min\mathfrak{D}$-estimator) if the
pseudodistance $\mathfrak{D}(P_{\theta},Q)$ is decomposable on $\mathcal{P}%
\otimes\mathcal{P}^{+}$ and the parameters $T_{\mathfrak{D}}(Q)\in\Theta$
minimize $\mathfrak{D}^{0}(P_{\theta})+\int\rho_{\theta}\mathrm{d}Q$ on
$\Theta,$ in symbols
\begin{equation}
T_{\mathfrak{D}}(Q)=\arg\inf_{\theta}\left[  \mathfrak{D}^{0}(P_{\theta}%
)+\int\rho_{\theta}\mathrm{d}Q\right] , \text{ \ \ \ for all \ }%
Q\in\mathcal{Q}.\label{c11}%
\end{equation}
In particular, for $Q=P_{n}\in\mathcal{P}_{\text{emp}}$
\begin{equation}
\widehat{\theta}_{\mathfrak{D},n}:=T_{\mathfrak{D}}(P_{n})=\arg\inf_{\theta
}\left[  \mathfrak{D}^{0}(P_{\theta})+\frac{1}{n}\sum_{i=1}^{n}\rho_{\theta
}(X_{i})\right] .\label{c1}%
\end{equation}
\end{definition}

\begin{theorem}
Every $\min\mathfrak{D}$-estimator
\begin{equation}
\widehat{\theta}_{\mathfrak{D},n}=\arg\inf_{\theta}\left[  \mathfrak{D}%
^{0}(P_{\theta})+\frac{1}{n}\sum_{i=1}^{n}\rho_{\theta}(X_{i})\right]
\label{c2}%
\end{equation}
is Fisher consistent in the sense that%
\begin{equation}
T_{\mathfrak{D}}(P_{\theta_{0}})=\theta_{0},\text{ \ \ for all }\theta_{0}%
\in\Theta. \label{c3}%
\end{equation}

\end{theorem}

\textit{Proof.} Consider arbitrary fixed $\theta_{0}\in\Theta$. Then, by
assumptions, $\mathfrak{D}^{1}(P_{\theta_{0}})$ is a finite constant.
Therefore (\ref{c11}) together with the definition of pseudodistance implies%
\begin{eqnarray*}
T_{\mathfrak{D}}(P_{\theta_{0}})  &  =&\arg\inf_{\theta}\left[  \mathfrak{D}%
^{0}(P_{\theta})+\int\rho_{\theta}\mathrm{d}P_{\theta_{0}}\right]
\\ &  =&\arg\inf_{\theta}\left[
\mathfrak{D}^{0}(P_{\theta})+\mathfrak{D}
^{1}(P_{\theta_{0}})+\int\rho_{\theta}\mathrm{d}P_{\theta_{0}}\right]
\\ &
=&\arg\inf_{\theta}\mathfrak{D}(P_{\theta},P_{\theta_{0}})=\theta_{0}.
\end{eqnarray*}

The decomposability of a pseudodistance $\mathfrak{D}(P_{\theta},Q)$ leads to
the additive structure of the empirical version
\begin{equation}
\mathfrak{D}(P_{\theta},P_{n})\sim\mathfrak{D}^{0}(P_{\theta})+\int
\rho_{\theta}\mathrm{d}P_{n}=\mathfrak{D}^{0}(P_{\theta})+\frac{1}{n}%
\sum_{i=1}^{n}\rho_{\theta}(X_{i})\label{c4}%
\end{equation}
in the definition (\ref{c2}) of the $\min\mathfrak{D}$-estimators,
which opens the possibility to apply the methods of the asymptotic
theory of $M$-estimators (cf. Hampel et al. (1986), van der Vaart
and Wellner (1996), van der Vaart (1998) or Mieske and Liese
(2008)).

\section{A special class of minimum pseudodistance estimators}

\subsection{Definitions and invariance properties}

For probability measures $P\in\mathcal{P}$ and $Q\in\mathcal{P}^{+}$\ consider
the following family of pseudodistances of orders $\alpha\geq0$,
\begin{equation}
\mathfrak{R}_{\alpha}(P,Q)=\frac{1}{1+\alpha}\ln\left(  \int p^{\alpha
}\mathrm{d}P\right)  +{\frac{1}{\alpha(1+\alpha)}}\ln\left(  \int q^{\alpha
}\mathrm{d}Q\right)  -{\frac{1}{\alpha}}\ln\left(  \int p^{\alpha}%
\mathrm{d}Q\right)  .\label{38}%
\end{equation}

The following basic condition which guarantees the finiteness of
the pseudodistances $\mathfrak{R}_{\alpha}(P,Q)$ is assumed. For
some positive $\beta,$
\begin{equation}
p^{\beta},\ q^{\beta},\ln p\in\mathbb{L}_{1}(Q)\text{ \ \ for all }%
P\in\mathcal{P},\text{ }Q\in\mathcal{P}^{+},\label{36B}%
\end{equation}
where $\mathbb{L}_{1}(Q):=\{f:\mathbb{R}^{d}\rightarrow\mathbb{R}\;\text{such
that}\;\int|f|\mathrm{d}Q<\infty\}$.

We then have:

\begin{theorem}
Let the condition (\ref{36B}) hold for some $\beta>0$. Then for all
$0<\alpha<\beta$, $\mathfrak{R}_{\alpha}(P,Q)$ defined in (\ref{38}) is a
family of pseudodistances decomposable in the sense
\begin{equation}
\mathfrak{R}_{\alpha}(P,Q)=\mathfrak{R}_{\alpha}^{0}(P)+\mathfrak{R}_{\alpha
}^{1}(Q)-{\frac{1}{\alpha}}\ln\left( \int p^{\alpha}\mathrm{d}Q\right)
,\label{38A}%
\end{equation}
where
\begin{equation}
\mathfrak{R}_{\alpha}^{0}(P)=\frac{1}{1+\alpha}\ln\left( \int p^{\alpha
}\mathrm{d}P\right) \text{ \ and \ }\mathfrak{R}_{\alpha}^{1}(Q)={\frac
{1}{\alpha(1+\alpha)}}\ln\left( \int q^{\alpha}\mathrm{d}Q\right) \label{38a}%
\end{equation}
and the limit relation holds
\begin{equation}
\mathfrak{R}_{\alpha}(P,Q)\rightarrow\mathfrak{R}_{0}(P,Q):=\int\ln
q\mathrm{d}Q-\int\ln p\mathrm{d}Q \text{ \ \ for }\alpha\downarrow
0.\label{38B}%
\end{equation}

\end{theorem}

\textit{Proof.} Under (\ref{36B}), the expressions $\ln(\int q^{\alpha
}\mathrm{d}Q),$\ $\ln(\int p^{\alpha}\mathrm{d}Q)$\ and $\int\ln p
\mathrm{d}Q$ appearing in (\ref{38}) and (\ref{38B}) are finite so that the
expressions $\mathfrak{R}_{\alpha}(P,Q)$ and $\mathfrak{R}_{0}(P,Q)$ are well
defined. Recall that, for arbitrary arguments $s,t>0$ and fixed parameters
$a,b>0$ with the property $1/a+1/b=1$ it holds
\begin{equation}
st\leq{\frac{s^{a}}{a}}+{\frac{t^{b}}{b}}\label{34a}%
\end{equation}
with equality if and only if $s^{a}=t^{b}$. Indeed, from the strict concavity
of the logarithmic function we deduce the inequality
\[
\ln(st)={\frac{1}{a}}\ln s^{a}+{\frac{1}{b}}\ln t^{b}\leq\ln\left(
{\frac{s^{a}}{a}}+{\frac{t^{b}}{b}}\right)
\]
and the stated condition for equality.

Taking $\alpha>0$ and substituting
\[
s=\frac{p^{\alpha}}{\left(  \int p^{\alpha a}\,\mathrm{d}\lambda\right)
^{1/a}},\text{ \ }t=\frac{q}{\left(  \int q^{b}\,\mathrm{d}\lambda\right)
^{1/b}}\text{ \ with \ }a=\frac{1+\alpha}{\alpha},\text{ \ }b=1+\alpha
\]
in the inequality (\ref{34a}), and integrating both sides by $\lambda$, we
obtain the inequality
\[
\int p^{\alpha}q\,\mathrm{d}\lambda\leq\left(  \int p^{1+\alpha}%
\,\mathrm{d}\lambda\right)  ^{\alpha/(1+\alpha)}\left(  \int q^{1+\alpha
}\,\mathrm{d}\lambda\right)  ^{1/(1+\alpha)}%
\]
with equality if and only if $p^{\alpha a}=q^{b}$ $\lambda$-a.s., i.e. if and
only if $p=q$ $\lambda$-a.s. Since the expression (\ref{38}) satisfies for
$\alpha>0$ the relation
\begin{equation}
\mathfrak{R}_{\alpha}(P,Q)=\frac{1}{\alpha}\left\{  \ln\left[
\left(  \int p^{1+\alpha}\,\mathrm{d}\lambda\right)
^{\alpha/(1+\alpha)}\left(  \int
q^{1+\alpha}\,\mathrm{d}\lambda\right)  ^{1/(1+\alpha)}\right]
-\ln\int p^{\alpha}q\,\mathrm{d}\lambda\right\},\label{dkd}
\end{equation}
we see that $\mathfrak{R}_{\alpha}(P,Q)$ is a pseudodistance on the space
$\mathcal{P}\otimes\mathcal{P}^{+}$. The decomposability in the sense of
(\ref{38A}) on this space is obvious and the limit relation
\[
\mathfrak{R}_{0}(P,Q)=\lim_{\alpha\downarrow0}\mathfrak{R}_{\alpha}(P,Q)
\]
results as follows:
\begin{eqnarray*}
&& \lim_{\alpha\downarrow0}\mathfrak{R}_{\alpha}(P,Q)=\\ &&
=\lim_{\alpha\downarrow0}\frac{1}{\alpha+1}\ln\left(  \int
p^{\alpha }\mathrm{d}P\right) +\frac{1}{\alpha(\alpha+1)}\ln\left(
\int q^{\alpha
}\mathrm{d}Q\right)  -\frac{1}{\alpha}\ln\left(  \int p^{\alpha}%
\mathrm{d}Q\right)  \\ &&
=\lim_{\alpha\downarrow0}\frac{1}{\alpha+1}\left[  \ln\left(  \int
p^{\alpha}\mathrm{d}P\right)  -\ln\left(  \int
q^{\alpha}\mathrm{d}P\right)
\right]  +\frac{1}{\alpha}\left[  \ln\left(  \int q^{\alpha}\mathrm{d}%
Q\right)  -\ln\left(  \int p^{\alpha}\mathrm{d}Q\right)  \right]  \\
&& =\lim_{\alpha\downarrow0}\frac{1}{\alpha}\ln\frac{\int q^{\alpha}%
\mathrm{d}Q}{\int p^{\alpha}\mathrm{d}Q}=\int\ln\frac{q}{p}\mathrm{d}%
Q=\mathfrak{R}_{0}(P,Q).
\end{eqnarray*}

Similarly as earlier in this paper, we are interested in the estimators
obtained by replacing the hypothetical probability measure $P_{\theta_{0}}$ in
the $\mathfrak{R}_{\alpha}$-pseudodistances $\mathfrak{R}_{\alpha}(P_{{\theta
}},P_{\theta_{0}})$ by the empirical measure $P_{n}$. Consider therefore the
family of minimum pseudodistance estimators of orders $0\leq\alpha\leq\beta$
(in symbols, $\min\mathfrak{R}_{\alpha}$-estimators) defined as $\widehat
{\theta}_{n}=T_{\alpha}(P_{n})$ for $T_{\alpha}(Q)\in\Theta$ with
$Q\in\mathcal{Q}=\mathcal{P}^{+}\cup\mathcal{P}_{\text{emp}}$ defined by
\begin{equation}
T_{\alpha}(Q)=\left\{
\begin{array}
[c]{ll}%
\arg\inf_{{\theta}}\left[  \frac{1}{1+\alpha}\ln\left(  \int p_{\theta
}^{\alpha}\mathrm{d}P_{\theta}\right)  -{\frac{1}{\alpha}}\ln(\int p_{\theta
}^{\alpha}\mathrm{d}Q)\right]   & \text{ \ \ \ }\mbox{if}\ 0<\alpha\leq
\beta\medskip\\
\arg\inf_{{\theta}}-\int\ln p_{{\theta}}\mathrm{d}Q & \text{ \ \ \ }%
\mbox{if}\text{ }\alpha=0.
\end{array}
\right.  \label{d1}%
\end{equation}
The upper formula is equivalent to
\begin{equation}
T_{\alpha}(Q)=\arg\sup_{\theta}\frac{\int p_{\theta}^{\alpha}\mathrm{d}%
Q}{C_{\alpha}(\theta)}\label{d3}%
\end{equation}
where
\begin{equation}
C_{\alpha}(\theta)=\left(  \int p_{\theta}^{1+\alpha}\mathrm{d}\lambda\right)
^{\alpha/(1+\alpha)}.\label{d2}%
\end{equation}

Hence, alternatively, we can write
\begin{equation}
\widehat{\theta}_{n}=\left\{
\begin{array}
[c]{ll}%
\arg\sup_{\theta}\left[ -\frac{1}{\alpha+1}\ln\left( \int p_{\theta}^{\alpha
}\mathrm{d}P_{\theta}\right) +\frac{1}{\alpha}\ln\left( \frac{1}{n}\sum
_{i=1}^{n}p_{\theta}^{\alpha}(X_{i})\right) \right]  & \text{ \ \ \ }%
\mbox{if}\ 0<\alpha\leq\beta\medskip\\
\arg\sup_{\theta}\frac{{ 1}}{{ n}}\sum_{i=1}^{n}\ln p_{{\theta}}(X_{i}) &
\text{ \ \ \ }\mbox{if}\text{ }\alpha=0
\end{array}
\right. \label{40r}%
\end{equation}
or
\begin{equation}
\widehat{\theta}_{n}=\left\{
\begin{array}
[c]{ll}%
\arg\sup_{\theta}C_{\alpha}{(\theta)}^{-1}\frac{{ 1}}{{ n}}\sum_{i=1}%
^{n}p_{{\theta}}^{\alpha}(X_{i}) & \text{ \ \ \ }\mbox{if}\ 0<\alpha\leq
\beta\medskip\\
\arg\sup_{\theta}\frac{{ 1}}{{ n}}\sum_{i=1}^{n}\ln p_{{\theta}}(X_{i}) &
\text{ \ \ \ }\mbox{if}\text{ }\alpha=0.
\end{array}
\right. \label{40}%
\end{equation}

Note that, for $\alpha\downarrow0$, the upper criterion function in
(\ref{40r}) tends to the lower maximum likelihood criterion. Indeed,
\begin{eqnarray*}
\lim_{\alpha\rightarrow0}\left[ -\frac{1}{\alpha+1}\ln\left( \int p_{\theta
}^{\alpha}\mathrm{d}P_{\theta}\right) +\frac{1}{\alpha}\ln\left( \frac{1}%
{n}\sum_{i=1}^{n}p_{\theta}^{\alpha}(X_{i})\right) \right] =\frac{{ 1}}{{ n}%
}\sum_{i=1}^{n}\ln p_{{\theta} }(X_{i})
\end{eqnarray*}
by the l'Hospital rule.

In the following, we give some invariance properties of $\min\mathfrak{R}%
_{\alpha}$-estimators.

If the statistical model $\langle(\mathbb{R}^{d},\mathcal{B}(\mathbb{R}%
^{d}));{\mathcal{P}}=(P_{\theta}:\theta\in\Theta)\rangle$ is reparametrized by
$\vartheta=\vartheta(\theta)$, then the new $\min\mathfrak{R}_{\alpha}%
$-estimators $\widehat{\vartheta}_{n}$ are related to the original
$\widehat{\theta}_{n}$ by $\widehat{\vartheta}_{n}=\vartheta(\widehat{\theta
}_{n})$. If the observations $x\in{\mathcal{X}}$ are replaced by $y=T(x)$,
where $T:(\mathbb{R}^{d},\mathcal{B}(\mathbb{R}^{d}))\mapsto(\mathbb{R}%
^{d},\mathcal{B}(\mathbb{R}^{d}))$ is a measurable statistic with the inverse
$T^{-1}$, then the densities
\[
\tilde{p}_{\theta}={\frac{\mathrm{d}\tilde{P}_{\theta}}{\mathrm{\ d}%
\tilde{\lambda}}}%
\]
in the transformed model $\tilde{{\mathcal{P}}}=(\tilde{P}_{\theta}=P_{\theta
}T^{-1}:\theta\in\Theta)$ with respect to $\sigma$-finite dominating measure
$\tilde{\lambda}=\lambda T^{-1}$ are related to the original densities
$p_{\theta}$ by
\begin{equation}
\tilde{p}_{\theta}(y)=p_{\theta}(T^{-1}y)\,{\mathcal{J}}_{T}(y) ,\label{c8}%
\end{equation}
where ${\mathcal{J}}_{T}(y)=\mathrm{d}\lambda T^{-1}/\mathrm{d}\tilde{\lambda
}$ is a generalized Jacobian of the statistic $T$. If $\lambda$\ is the
Lebesque measure and the inverse mapping $H=T^{-1}$ is differentiable, then
${\mathcal{J}}_{T}(y)$ is the determinant
\[
{\mathcal{J}}_{T}(y)=\left\vert \frac{\mathrm{d}}{\mathrm{d}y}H(y)\right\vert
.
\]

The $\min\mathfrak{R}_{\alpha}$-estimators are in general not equivariant with
respect to invertible transformations $T$ of observations, unless $\alpha=0$.

\begin{theorem}
The $\min\mathfrak{R}_{\alpha}$-estimators $\tilde{\theta}_{n}$ in the above
considered transformed model coincide with the original $\min\mathfrak{R}%
_{\alpha}$-estimators $\widehat{\theta}_{n}$, if the Jacobian ${\mathcal{J}%
}_{T}$ of transformation is a nonzero constant on the transformed observation
space. Thus, the $\min\mathfrak{R}_{\alpha}$-estimators are equivariant under
linear statistics $Tx=ax+b$.
\end{theorem}

\textit{Proof.} For $\alpha=0$ the $\min\mathfrak{R}_{\alpha}$-estimator
coincides with the maximum likelihood estimator, whose equivariance is well
known. For $\alpha>0$, by (\ref{c8}) and (\ref{40}),
\begin{eqnarray*}
\tilde{\theta}_{n}  & =&\arg\sup_{\theta}C_{\alpha}{(\theta)}^{-1}\frac{{1}%
}{{n}}\sum_{i=1}^{n}\tilde{p}_{\theta}^{\alpha}(TX_{i})\\
& =&\arg\sup_{\theta}C_{\alpha}{(\theta)}^{-1}\frac{{1}}{{n}}\sum_{i=1}%
^{n}p_{\theta}^{\alpha}(X_{i})\,{\mathcal{J}}_{T}^{\alpha}(TX_{i}).
\end{eqnarray*}
Comparing with (\ref{40}) it follows that $\tilde{\theta}_{n}=\widehat{\theta
}_{n}$ if $y\mapsto{\mathcal{J}}_{T}(y)$ is a nonzero constant. If $\alpha=0$,
then the estimator coincides with the maximum likelihood estimator and its
equivariance is well known.

\subsection{Limit properties of $\min\mathfrak{R}_{\alpha}$-estimators}

Define
\[
R_{\alpha}(\theta_{0}):=\sup_{\theta}\int h(x,\theta)\mathrm{d}P_{\theta_{0}%
}(x)
\]
where
\[
h(x,\theta):=\frac{p_{\theta}^{\alpha}(x)}{C_{\alpha}\left(  \theta\right)
}.
\]

By the Fisher consistency of the functional $T_{\alpha}$ defined in
(\ref{d3}), it holds
\[
\arg\sup_{\theta}\int h(x,\theta)\mathrm{d}P_{\theta_{0}}(x)=\theta_{0}%
\]
and $\theta_{0}$ is the only optimizer in the above expression, since
$\mathfrak{R}_{\alpha}\left(  P_{\theta},P_{\theta_{0}}\right)  =0$ implies
$\theta=\theta_{0}.$

Define the estimate of $R_{\alpha}(\theta_{0})$ through
\[
\widehat{R}_{\alpha}(\theta_{0}):=\sup_{\theta}\int h(x,\theta)\mathrm{d}%
P_{n}=\sup_{\theta}\frac{1}{n}\sum_{i=1}^{n}h(X_{i},\theta),
\]
where the $\theta_{0}$ indicates that the sampling is i.i.d. under
$P_{\theta_{0}}.$ The $\min\mathfrak{R}_{\alpha}$-estimator $\widehat{\theta
}_{n}$ is then defined through
\[
\widehat{\theta}_{n}=\arg\sup_{\theta}\frac{1}{n}\sum_{i=1}^{n}h(X_{i}%
,\theta).
\]
This optimum need not be uniquely defined.

The usual regularity properties of the model will be assumed
throughout the rest of the paper, namely:  (i) The density
$p_{\theta}(x)$ has continuous partial derivatives with respect to
$\theta$ up to 3th order (for all $x$ $\lambda$-a.e.). (ii) There
exists a neighborhood $N(\theta_{0})$ of $\theta_{0}$ such that
the first, the second and the third order partial derivatives
(w.r.t. $\theta$) of $h(x,\theta)$ are dominated on $N(\theta
_{0})$ by some $P_{\theta_{0}}$-integrable functions. (iii) The
integrals
$\int(\partial^{2}/\partial\theta^{2})h(x,\theta_{0})\mathrm{d}P_{\theta_{0}%
}(x)$ and $\int(\partial/\partial\theta)h(x,\theta_{0})(\partial
/\partial\theta)^{t}h(x,\theta_{0})\mathrm{d}P_{\theta_{0}}(x)$ exist.

\begin{theorem}
\label{thasympt} Assume that the above conditions hold.

\begin{enumerate}
\item[(a)] Let $B:=\left\{  \theta\in\Theta;~\Vert\theta-\theta_{0}\Vert\leq
n^{-1/3}\right\}  $. Then, as $n\rightarrow\infty$, with probability one, the
function $\theta\mapsto\frac{1}{n}\sum_{i=1}^{n}h(X_{i},\theta)$ attains a
local maximal value at some point $\widehat{\theta}_{n}$ in the interior of
$B$, which implies that the estimate $\widehat{\theta}_{n}$ is $n^{1/3}$-consistent.

\item[(b)] $\sqrt{n}\left(  \widehat{\theta}_{n}-\theta_{0}\right)  $
converges in distribution to a centered multivariate normal random variable
with covariance matrix
\begin{equation}
V=S^{-1}MS^{-1}\label{variance limite 1}%
\end{equation}
with $S:=-\int(\partial^{2}/\partial\theta^{2})h(x,\theta_{0})\mathrm{d}%
P_{\theta_{0}}(x)$ and $M:=\int(\partial/\partial\theta)h(x,\theta
_{0})(\partial/\partial\theta)^{t}h(x,\theta_{0})\mathrm{d}P_{\theta_{0}}(x)$.

\item[(c)] $\sqrt{n}\left(  \widehat{R}_{\alpha}(\theta_{0})-R_{\alpha}%
(\theta_{0})\right)  $ converges in distribution to a centered normal variable
with variance $\sigma^{2}(\theta_{0})=\int h(x,\theta_{0})^{2}\mathrm{d}%
P_{\theta_{0}}(x)-\left(  \int h(x,\theta_{0})\mathrm{d}P_{\theta_{0}%
}(x)\right)  ^{2}$.
\end{enumerate}
\end{theorem}

\textit{Proof.} (a) A simple calculus give
\begin{equation}
\int(\partial/\partial\theta)h(x,\theta_{0})\mathrm{d}P_{\theta_{0}%
}(x)=0\label{eqn1}%
\end{equation}
and
\begin{equation}
\int(\partial^{2}/\partial\theta^{2})h(x,\theta_{0})\mathrm{d}P_{\theta_{0}%
}(x)=-S.\label{eqn2}%
\end{equation}
Observe that the matrix $S$ is symmetric and positive definite.

Let $U_{n}:=\frac{1}{n}\sum_{i=1}^{n}(\partial/\partial\theta) h(X_{i}%
,\theta_{0})$, and use (\ref{eqn1}) in connection with the central limit
theorem to see that
\begin{equation}
\sqrt{n}U_{n}\rightarrow\mathcal{N}(0,M).\label{eqn3}%
\end{equation}
Also, let $V_{n}:=\frac{1}{n}\sum_{i=1}^{n}(\partial^{2}/\partial\theta
^{2})h(X_{i},\theta_{0})$, and use (\ref{eqn2}) in connection with the law of
large numbers to conclude that
\begin{equation}
V_{n}\rightarrow-S~~(a.s).\label{eqn4}%
\end{equation}
Now, for any $\theta=\theta_{0}+un^{-1/3}$ with $|u|\leq1$, consider a Taylor
expansion of $\frac{1}{n}\sum_{i=1}^{n}h(X_{i},\theta)$ in $\theta$ around
$\theta_{0}$, and use the hypothesis to see that
\[
\sum_{i=1}^{n}h(X_{i},\theta)-\sum_{i=1}^{n}h(X_{i},\theta_{0})=n^{2/3}%
u^{t}U_{n}+2^{-1}n^{1/3}u^{t}V_{n}u+O(1)~(a.s.)
\]
uniformly on $u$ with $|u|\leq1$. Now, use (\ref{eqn4}) and the fact that
$U_{n}=O\left(  n^{-1/2}(\log\log n)^{1/2}\right)  $ (a.s) to conclude that
\[
\sum_{i=1}^{n}h(X_{i},\theta)-\sum_{i=1}^{n}h(X_{i},\theta_{0})=O\left(
n^{1/6}(\log\log n)^{1/2}\right)  -2^{-1}u^{t}Sun^{1/3}+O(1)~(a.s.)
\]
uniformly on $u$ with $|u|\leq1$. Hence, uniformly on the surface of the ball
$B$ (i.e., uniformly on $u$ with $|u|=1$), we have
\begin{equation}
nP_{n}h(\theta)-nP_{n}h(\theta_{0})\leq O\left(  n^{1/6}(\log\log
n)^{1/2}\right)  -2^{-1}cn^{1/3}+O(1)~~(a.s.)\label{eqn5}%
\end{equation}
where $c$ is the smallest eigenvalue of the matrix $S$. Note that $c$ is
positive since $S$ is positive definite (it is symmetric, positive and non
singular). In view of $(\ref{eqn5})$, by the continuity of $\theta\mapsto
\sum_{i=1}^{n}h(X_{i},\theta)-\sum_{i=1}^{n}h(X_{i},\theta_{0})$ and since it
takes value zero on $\theta=\theta_{0}$ and is asymptotically nonpositive, it
holds that as $n\rightarrow\infty$, with probability one, $\theta\mapsto
\frac{1}{n}\sum_{i=1}^{n}h(X_{i},\theta)$ attains its maximum value at some
point $\widehat{\theta}_{n}$ in the interior of the ball $B$, and therefore
the estimate $\widehat{\theta}_{n}$ satisfies $\left\Vert \widehat{\theta}%
_{n}-\theta_{0}\right\Vert =O(n^{-1/3})$.

(b) Using the fact that $\frac{1}{n}\sum_{i=1}^{n}(\partial/\partial
\theta)h(X_{i},\widehat{\theta}_{n})=0$ and a Taylor expansion of $\frac{1}%
{n}\sum_{i=1}^{n}(\partial/\partial\theta)h(X_{i},\theta)$ in $\widehat
{\theta}_{n}$ around $\theta_{0}$, we obtain
\[
0=\frac{1}{n}\sum_{i=1}^{n}(\partial/\partial\theta)h(X_{i},\widehat{\theta
}_{n})=\frac{1}{n}\sum_{i=1}^{n}(\partial/\partial\theta)h(X_{i},\theta
_{0})+(\widehat{\theta}_{n}-\theta_{0})^{t}\sum_{i=1}^{n} (\partial
^{2}/\partial\theta^{2})h(X_{i},\theta_{0})+o_{p}(n^{-1/2}).
\]
Hence,
\begin{equation}
\sqrt{n}\left(  \widehat{\theta}_{n}-\theta_{0}\right)  =-V_{n}^{-1}\sqrt
{n}U_{n}+o_{p}(1).\label{eqn6}%
\end{equation}
Using (\ref{eqn3}) and (\ref{eqn4}) and Slutsky theorem, we conclude then
\[
\sqrt{n}\left(  \widehat{\theta}-\theta_{0}\right)  \rightarrow\mathcal{N}%
\left(  0,V\right)
\]
where $V=S^{-1}M S^{-1}$.

(c) A Taylor expansion of $\widehat{R}_{\alpha}\left(  \theta_{0}\right)
=\frac{1}{n}\sum_{i=1}^{n}h(X_{i},\theta)$ in $\widehat{\theta}_{n}$ around
$\theta_{0}$, using the fact that $\int(\partial/\partial\theta)h(x,\theta
_{0})\mathrm{d}P_{\theta_{0}}(x)=0$, gives
\begin{equation}
\widehat{R}_{\alpha}\left(  \theta_{0}\right)  =\frac{1}{n}\sum_{i=1}%
^{n}h(X_{i},\theta_{0})+o_{p}(n^{-1/2}).
\end{equation}
Hence,
\[
\sqrt{n}\left(  \widehat{R}_{\alpha}\left(  \theta_{0}\right)  -R_{\alpha
}(\theta_{0})\right)  =\sqrt{n}\left[  \frac{1}{n}\sum_{i=1}^{n}h(X_{i}%
,\theta_{0})-\int h(x,\theta_{0})\mathrm{d}P_{\theta_{0}}(x)\right]
+o_{p}(1),
\]
which by the central limit theorem, converges in distribution to a centered
normal variable with variance $\sigma^{2}(\theta_{0})=\int h(x,\theta_{0}%
)^{2}\mathrm{d}P_{\theta_{0}}(x)-\left(  \int h(x,\theta_{0})\mathrm{d}%
P_{\theta_{0}}(x)\right)  ^{2}$.

\subsection{Robustness results}

We recall that a map $T$ defined on a set of probability measures and the
parameter space valued is a statistical functional corresponding to an
estimator $\widehat{\theta}_{n}$ of the parameter $\theta_{0}$, whenever
$T(P_{n})=\widehat{\theta}_{n}$. The influence function of the functional $T$
in $P$ measures the effect on $T$ of adding a small mass at $x$ and is defined
as
\[
\mathrm{IF}(x;T,P)=\lim_{\varepsilon\rightarrow0}\frac{T(\widetilde
{P}_{\varepsilon x})-T(P)}{\varepsilon},
\]
where $\widetilde{P}_{\varepsilon x}=(1-\varepsilon)P+\varepsilon\delta_{x}$.
When the influence function is bounded, the corresponding estimator is called
B-robust. The gross error sensitivity of $T$ at $P$ is defined by
\begin{equation}
\gamma^{*}(T,P):=\sup_{x}\|\mathrm{IF}(x;T,P)\|,
\end{equation}
the supremum being taken over all $x$ where $\mathrm{IF}(x;T,P)$ exists.

The statistical functional corresponding to the $\min\mathfrak{R}_{\alpha}%
$-estimator is
\begin{equation}
T_{\alpha}(Q):=\arg\sup_{\theta}\frac{\int p_{\theta}^{\alpha}\mathrm{d}%
Q}{C_{\alpha}(\theta)},
\end{equation}
where $C_{\alpha}(\theta):=(\int p_{\theta}^{\alpha+1}d\lambda)^{\frac{\alpha
}{\alpha+1}}$. Derivation w.r.t. $\theta$ shows that $T_{\alpha}(Q)$ is a
solution of the equation
\begin{equation}
\int[p_{\theta}^{\alpha-1}\dot{p}_{\theta}-c_{\alpha}(\theta)p_{\theta
}^{\alpha}]\mathrm{d}Q=0,
\end{equation}
where $c_{\alpha}(\theta):=\frac{\int p_{\theta}^{\alpha}\dot{p}_{\theta
}d\lambda}{\int p_{\theta}^{\alpha+1}d\lambda}$ and $\dot{p}_{\theta}$ is the
derivative with respect to $\theta$ of $p_{\theta}$.

\begin{theorem}
The influence function of $T_{\alpha}$ in $P_{\theta}$ is given by
\begin{equation}
\label{IF}\mathrm{IF}(x;T_{\alpha},P_{\theta})=M_{\alpha}(\theta
)^{-1}[p_{\theta}^{\alpha-1}(x)\dot{p}_{\theta}(x)-c_{\alpha}(\theta
)p_{\theta}^{\alpha}(x)],
\end{equation}
where
\begin{equation}
M_{\alpha}(\theta)=\int p_{\theta}^{\alpha-1}\dot{p}_{\theta}\dot{p}_{\theta
}^{t}d \lambda-\frac{\int p_{\theta}^{\alpha}\dot{p}_{\theta}d\lambda(\int
p_{\theta}^{\alpha}\dot{p}_{\theta}d\lambda)^{t} }{\int p_{\theta}^{\alpha
+1}d\lambda}.\label{Malpha}%
\end{equation}

\end{theorem}

\textit{Proof.} Consider the contaminated model $\widetilde{P}_{\varepsilon
x}=(1-\varepsilon)P_{\theta}-\varepsilon\delta_{x}$. Then
\begin{eqnarray*}
&& (1-\varepsilon)\int\{p_{T_{\alpha}(\widetilde{P}_{\varepsilon x})}%
^{\alpha-1}\dot{p}_{T_{\alpha}(\widetilde{P}_{\varepsilon
x})}-c_{\alpha }(T_{\alpha}(\widetilde{P}_{\varepsilon
x}))p_{T_{\alpha}(\widetilde {P}_{\varepsilon
x})}^{\alpha}\}\mathrm{d}P_{\theta}+\\ &&
+\varepsilon\{p_{T_{\alpha}(\widetilde{P}_{\varepsilon
x})}^{\alpha -1}(x)\dot{p}_{T_{\alpha}(\widetilde{P}_{\varepsilon
x})}(x)-c_{\alpha }(T_{\alpha}(\widetilde{P}_{\varepsilon
x}))p_{T_{\alpha}(\widetilde {P}_{\varepsilon
x})}^{\alpha}(x)\}=0.
\end{eqnarray*}
Derivating with respect to $\varepsilon$ in the above display and taking the
derivative in $\varepsilon=0$, after some calculations, we obtain
\begin{equation}
\mathrm{IF}(x;T_{\alpha},P_{\theta})=M_{\alpha}(\theta)^{-1}[p_{\theta
}^{\alpha-1}(x)\dot{p}_{\theta}(x)-c_{\alpha}(\theta)p_{\theta}^{\alpha}(x)],
\end{equation}
where $M_{\alpha}(\theta)$ is given by the formula (\ref{Malpha}).

Beside the influence function, the breakdown point provides information about
the robustness of an estimator. The breakdown point of an estimator
$\widehat{\theta}_{n}$ of a parameter $\theta_{0}$ is the largest amount of
contamination that the data may contain, such that $\widehat{\theta}_{n}$
still gives some information about $\theta_{0}$. Following Maronna et al.
(2006) (p.58), the asymptotic contamination breakdown point of an estimator
$\widehat{\theta}_{n}$ at $P_{\theta_{0}}$, denoted by $\varepsilon
^{*}(\widehat{\theta}_{n},\theta_{0})$, is the largest $\varepsilon^{*}%
\in(0,1)$ such that for $\varepsilon<\varepsilon^{*}$, $T((1-\varepsilon
)P_{\theta_{0}}+\varepsilon P)$ as function of $P$ remains bounded and also
bounded away from the boundary $\partial\Theta$ of $\Theta$. Here,
$T((1-\varepsilon)P_{\theta_{0}}+\varepsilon P)$ is the asymptotic value of
the estimator at $(1-\varepsilon)P_{\theta_{0}}+\varepsilon P$ by means of the
convergence in probability. The definition means that there exists a bounded
and closed set $K\subset\Theta$ such that $K\cap\partial\Theta=\emptyset$ and
\[
T((1-\varepsilon)P_{\theta_{0}}+\varepsilon P)\in K,\;\text{for all}%
\;\varepsilon<\varepsilon^{*}\;\text{and all}\; P.
\]

\subsubsection{Scale models}

\begin{center}
{\scriptsize \begin{figure}[ptb]
{\scriptsize
\begin{tabular}
[c]{cc}%
\includegraphics[width=7.5cm,height=7.5cm]{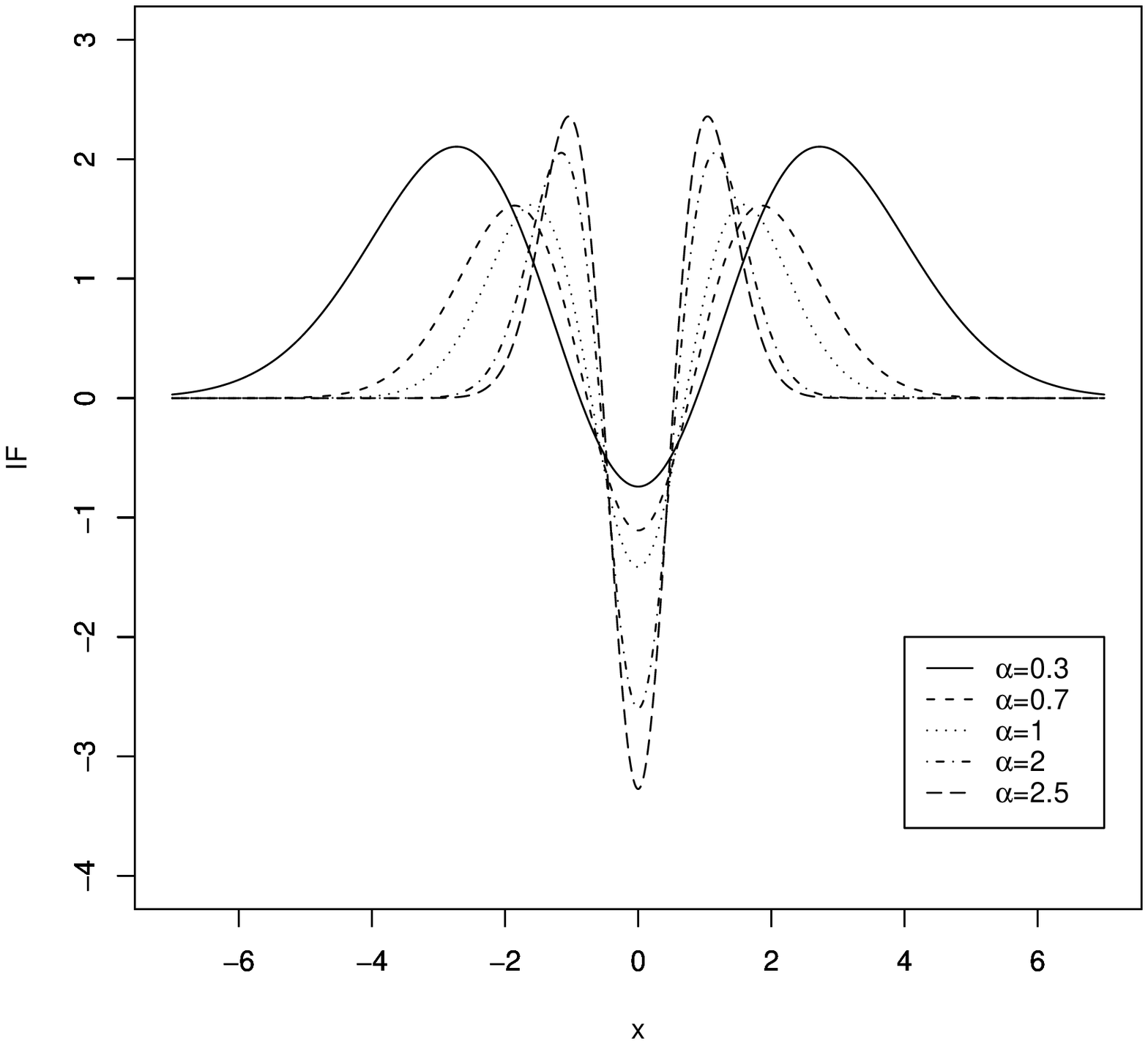} &
\includegraphics[width=7.5cm,height=7.5cm]{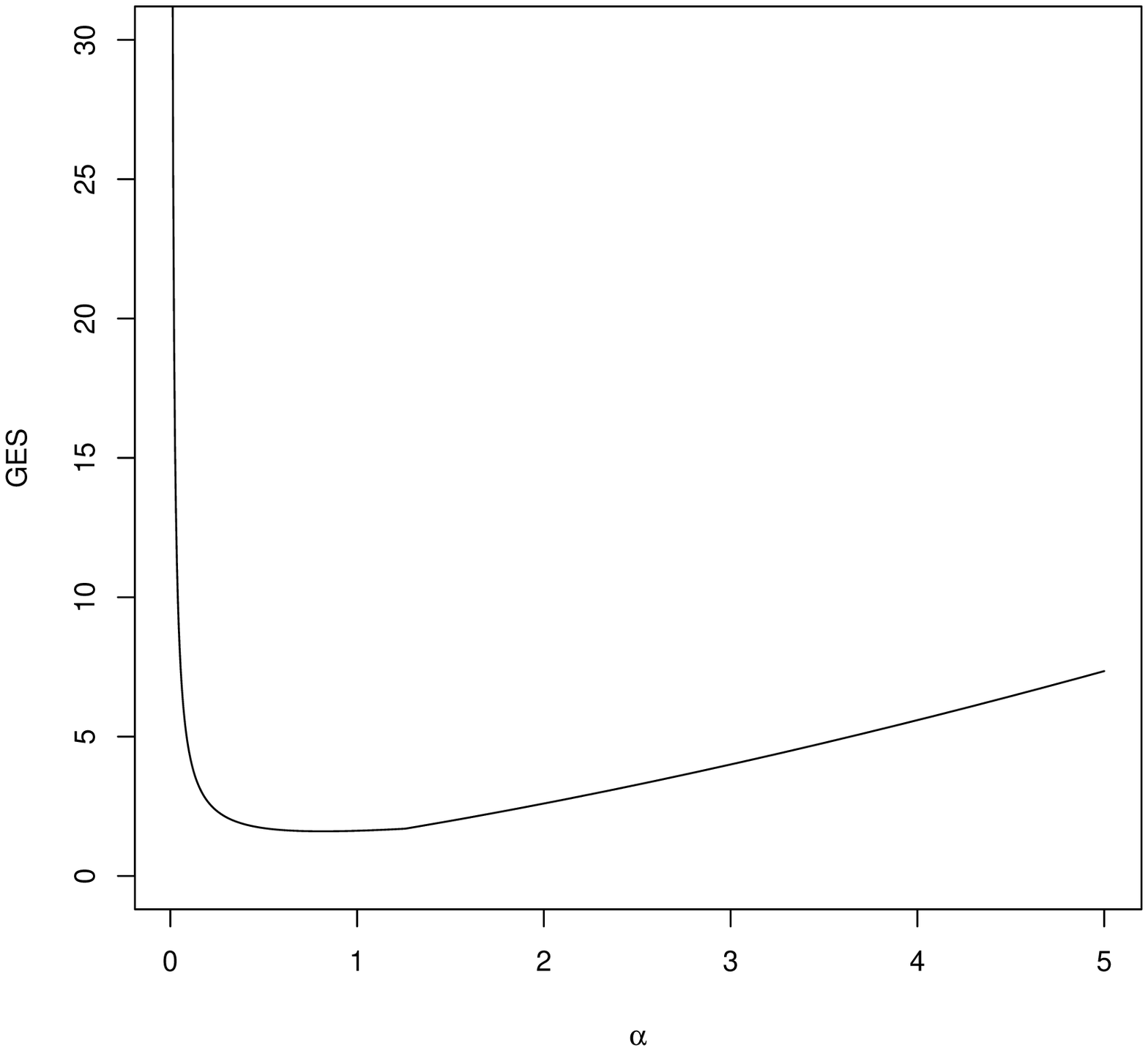}
\end{tabular}
\vspace{-2mm}}\caption{Influence functions and gross error sensitivity of
$\min\mathfrak{R}_{\alpha}$-estimators of the scale parameter $\sigma=1$ from
the normal model, when $m=0$ is known}%
\label{Gr1}%
\end{figure}}
\end{center}

\textit{(a) Standard deviation of univariate normal.} Consider the scale
normal model with known mean $m$ and take $\theta=\sigma$. The influence
function (\ref{IF}) takes on the form
\begin{equation}
\mathrm{IF}(x;T_{\alpha},P_{\sigma})=\frac{\sigma(\alpha+1)^{5/2}}{2}\left[
\left(  \frac{x-m}{\sigma}\right)  ^{2}-\frac{1}{\alpha+1}\right]  \exp\left(
-\frac{\alpha}{2}\left(  \frac{x-m}{\sigma}\right)  ^{2}\right)  .
\end{equation}

The gross error sensitivity of $T_{\alpha}$ in $P_{\sigma}$ is given by
\begin{equation}
\gamma^{*}(T_{\alpha},P_{\sigma})=\max\left\{  \frac{\sigma(\alpha+1)^{3/2}%
}{2},\frac{\sigma(\alpha+1)^{5/2}}{\alpha}\exp\left(  -\frac{3\alpha
+2}{2(\alpha+1)}\right)  \right\}  ,
\end{equation}
independently upon the value of $m$.

\begin{center}
{\scriptsize \begin{figure}[ptb]
{\scriptsize
\begin{tabular}
[c]{cc}%
\includegraphics[width=7.5cm,height=7.5cm]{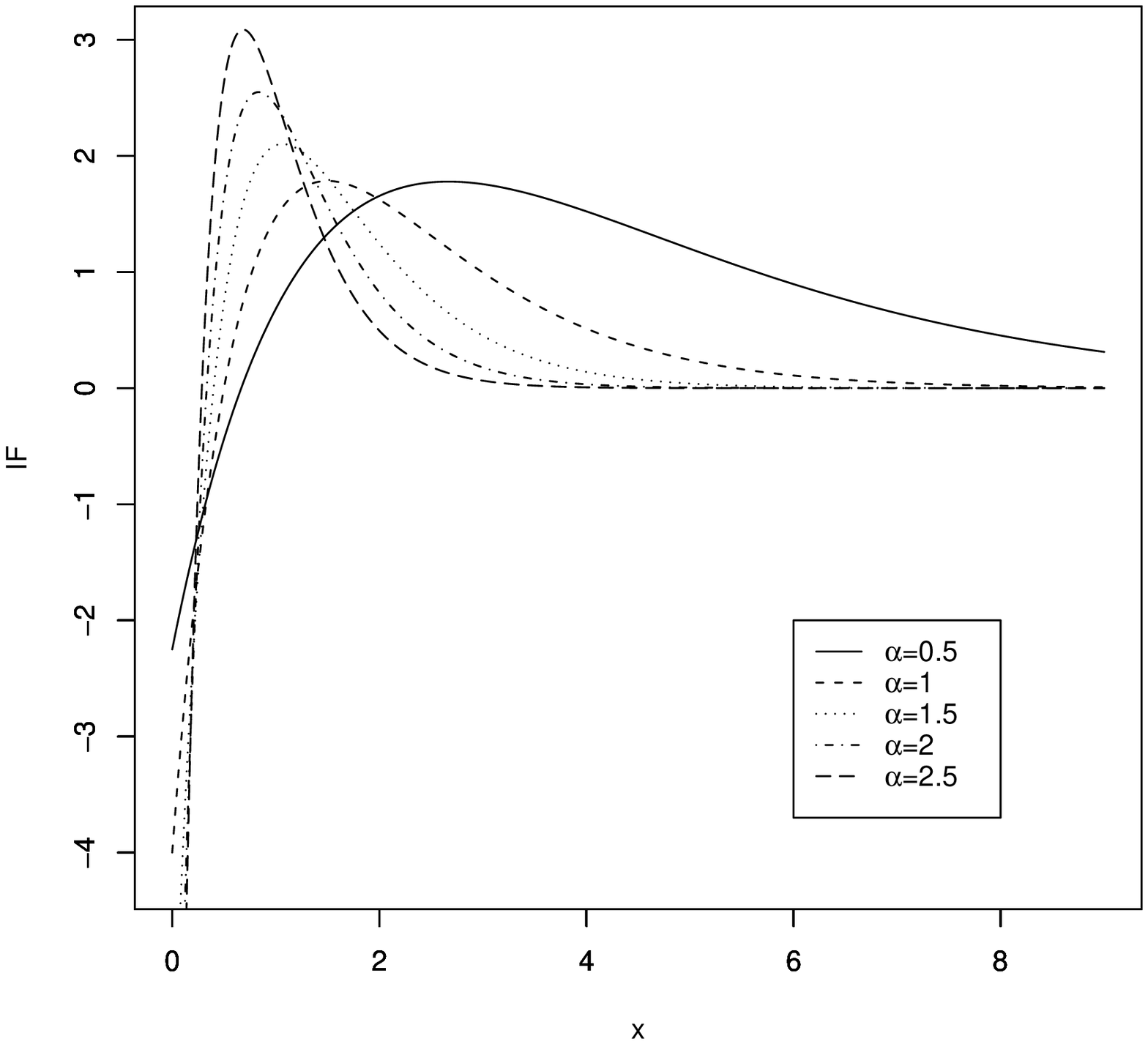} &
\includegraphics[width=7.5cm,height=7.5cm]{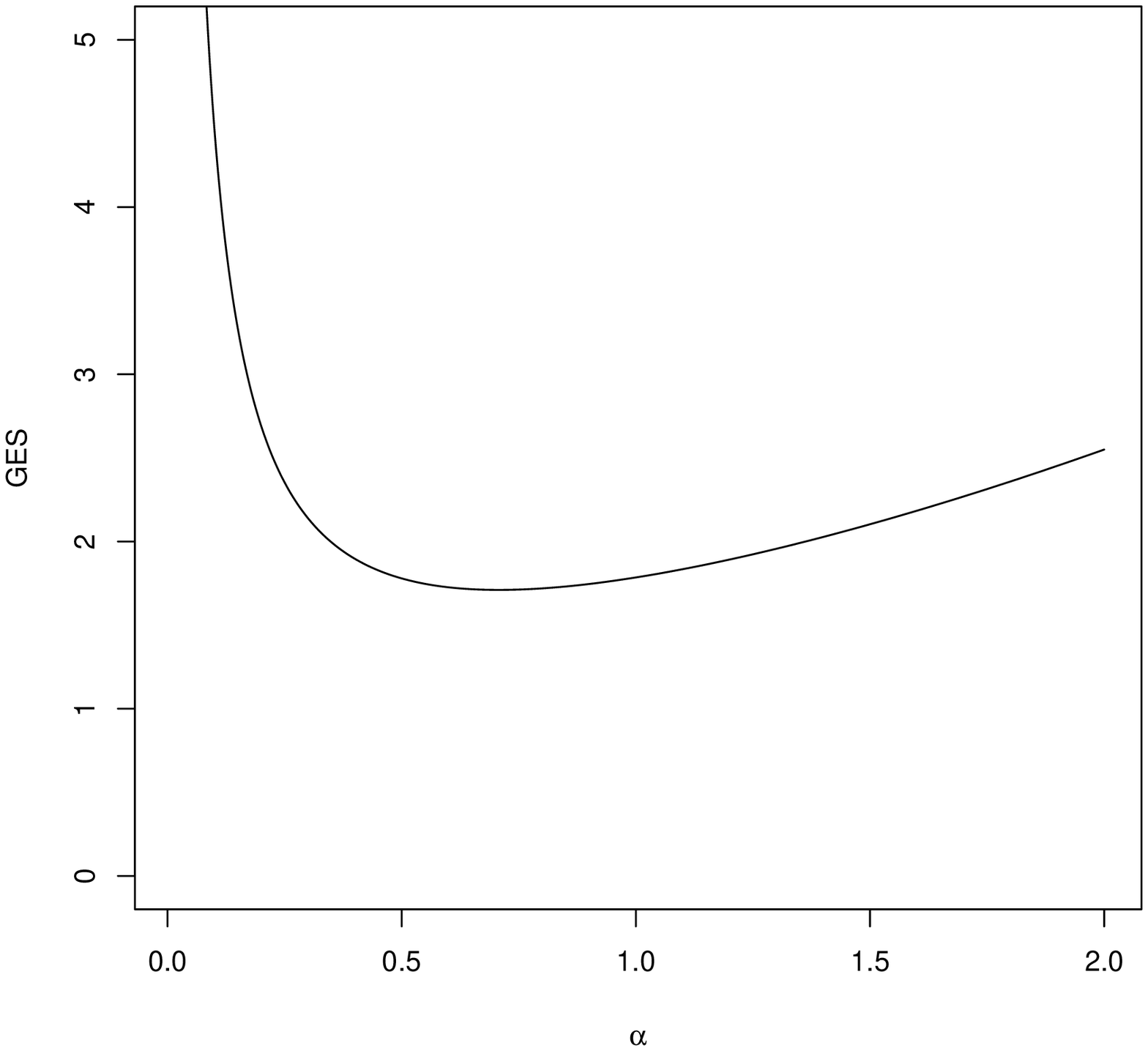}
\end{tabular}
\vspace{-2mm} }\caption{Influence functions and gross error sensitivity of
$\min\mathfrak{R}_{\alpha}$-estimators of the scale parameter $\theta=1$ from
the exponential model}%
\label{Gr2}%
\end{figure}}
\end{center}

Figure \ref{Gr1} presents influence functions $\mathrm{IF}(x;T_{\alpha
},P_{\sigma})$, in the case of the scale normal model with the known mean
$m=0$, when $\sigma=1$. For the same model, the gross error sensitivity of the
$\min\mathfrak{R}_{\alpha}$-estimator, as function of $\alpha$, is
represented. The gross error sensitivity attains its minimum value $\gamma
^{*}(T_{\alpha},P_{\sigma})=1.600413$, for $\alpha=0.81648$. This means that
the $\min\mathfrak{R}_{\alpha}$-estimator corresponding to $\alpha=0.81648$ is
the most B-robust estimator within the class of the $\min\mathfrak{R}_{\alpha
}$-estimators of $\sigma=1$.

On the other hand, the asymptotic relative efficiency of $T_{\alpha}$ is
\begin{equation}
\mathrm{ARE}(T_{\alpha},P_{\sigma})=\frac{2(2\alpha+1)^{5/2}}{(\alpha
+1)^{3}(3\alpha^{2}+4\alpha+2)}.
\end{equation}
Results for different values of $\alpha$ are given in the first row of Table
1. Note that, when $\alpha$ increases, the asymptotic relative efficiency of
the estimator decreases. Therefore, positive values of $\alpha$ close to zero
will ensure high efficiency and in the meantime the B-robustness of the
estimator. For example, the asymptotic relative efficiency of the estimator is
0.975 for $\alpha=0.1$, respectively 0.919 for $\alpha=0.2$. For both cases,
the $\min\mathfrak{R}_{\alpha}$-estimator is B-robust.

\textit{(b) Exponential model.} Consider the exponential model with density
$p_{\theta}(x)=\frac{1}{\theta}\exp(-\frac{x}{\theta})$, $x\geq0$. The
influence function of a $\min\mathfrak{R}_{\alpha}$-estimator of the parameter
$\theta$ is
\begin{equation}
\mathrm{IF}(x;T_{\alpha},P_{\theta})=\theta(\alpha+1)^{3}\left(  \frac
{x}{\theta}-\frac{1}{\alpha+1}\right)  \exp\left(  -\frac{\alpha x }{\theta
}\right)
\end{equation}
and the corresponding gross error sensitivity is
\begin{equation}
\label{GESexp}\gamma^{*}(T_{\alpha},P_{\theta})=\theta\frac{(\alpha+1)^{3}%
}{\alpha}\exp\left(  -\frac{2\alpha+1}{\alpha+1}\right)  .
\end{equation}

In Figure \ref{Gr2}, for different values of $\alpha$, influence functions of
$\min\mathfrak{R}_{\alpha}$-estimators of the parameter $\theta=1$ from the
exponential model are represented. The gross error sensitivity (\ref{GESexp})
as function of $\alpha$ is also represented. The most B-robust estimator over
the class of $\min\mathfrak{R}_{\alpha}$-estimators of $\theta=1$ is
associated to $\alpha=0.707$, case in which the gross error sensitivity takes
on the value $1.710$.

The asymptotic relative efficiency of $T_{\alpha}$ in $P_{\theta}$ is given
by
\begin{equation}
\mathrm{ARE}(T_{\alpha},P_{\theta})=\frac{(2\alpha+1)^{3}}{(\alpha
+1)^{4}(2\alpha^{2}+2\alpha+1)}.
\end{equation}
As in the case of the scale normal model, the efficiency remains high for
small $\alpha$. Thus, positive values of $\alpha$ close to zero will ensure
high efficiency and the B-robustness of the estimation procedure.

\subsubsection{Location models}

\begin{center}
{\scriptsize \begin{figure}[ptb]
{\scriptsize
\begin{tabular}
[c]{cc}%
\includegraphics[width=7.5cm,height=7.5cm]{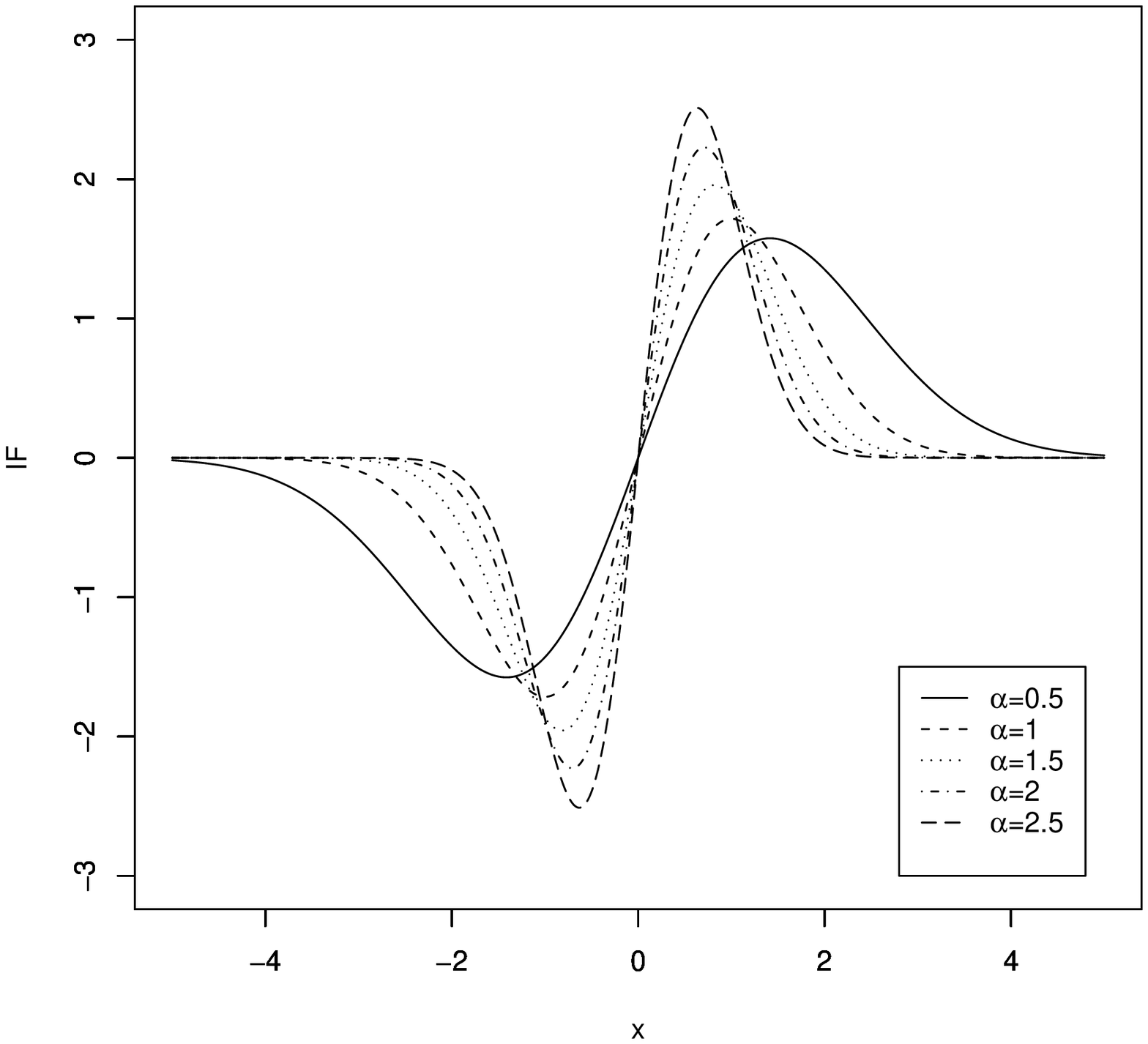} &
\includegraphics[width=7.5cm,height=7.5cm]{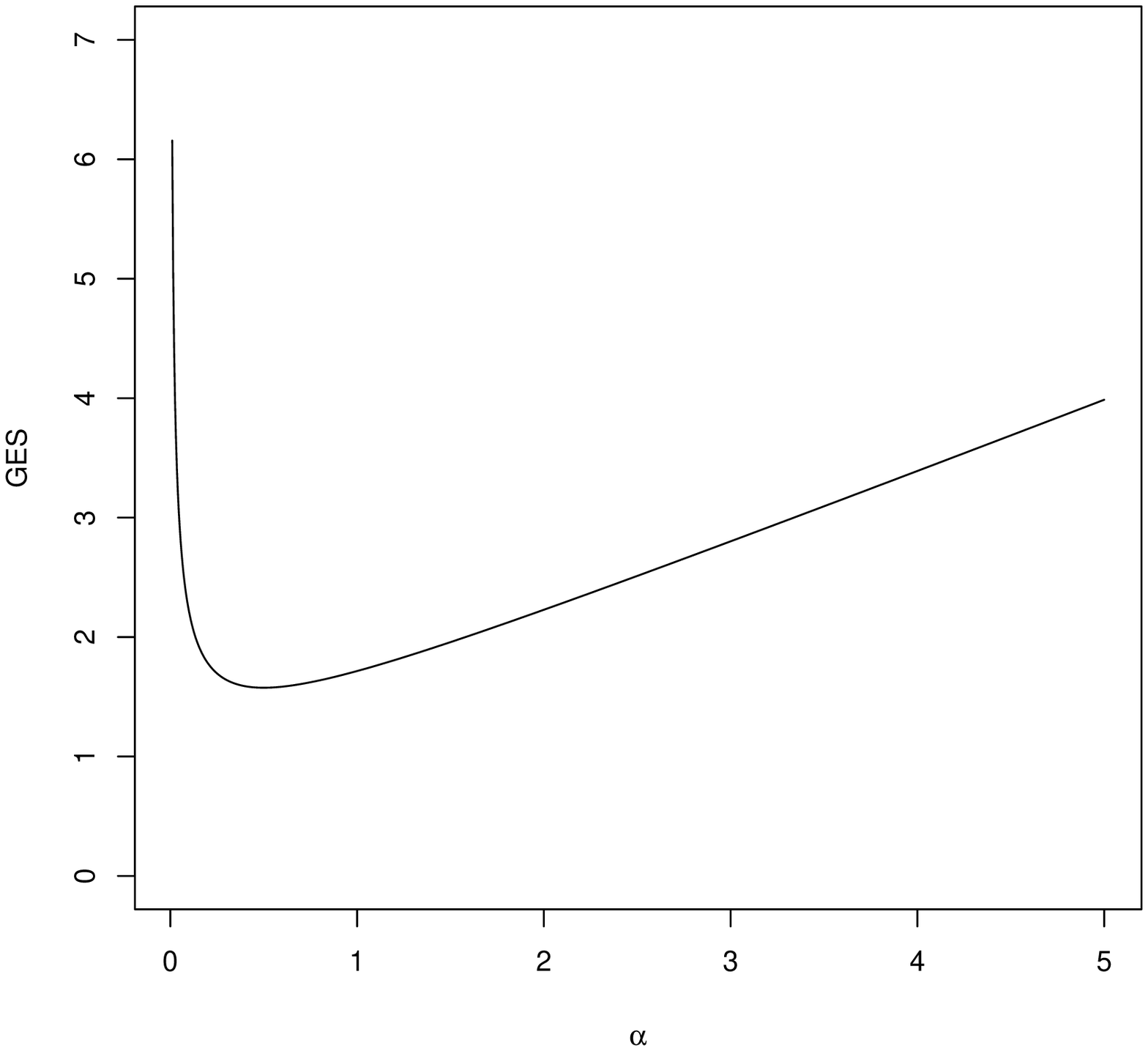}
\end{tabular}
\vspace{-2mm} }\caption{Influence functions and gross error sensitivity of
$\min\mathfrak{R}_{\alpha}$-estimators of the location parameter $m=0$ from
the normal model, when $\sigma=1$ is known}%
\label{Gr3}%
\end{figure}}
\end{center}

\textit{(c) Mean of univariate normal}. Letting $p_{\theta}$ be the
$\mathcal{N}(m,\sigma)$ density with known $\sigma$, the influence function of
a $\min\mathfrak{R}_{\alpha}$-estimator of the location parameter $\theta=m$
is
\begin{equation}
\mathrm{IF}(x;T_{\alpha},P_{m})=(\alpha+1)^{3/2}(x-m)\exp\left(  -\frac
{\alpha}{2}\left(  \frac{x-m}{\sigma}\right)  ^{2}\right)
\end{equation}
and the gross error sensitivity is
\begin{equation}
\label{GESlocnor}\gamma^{*}(T_{\alpha},P_{m})=(\alpha+1)^{3/2}\frac{\sigma
}{\sqrt{\alpha}}\exp(-1/2),
\end{equation}
independently from the value of $m$.

In Figure \ref{Gr3}, for $\sigma=1$ and different values of $\alpha$,
influence functions of $\min\mathfrak{R}_{\alpha}$-estimators of the location
parameter $m=0$ are represented. The gross error sensitivity (\ref{GESlocnor})
as function of $\alpha$ is also represented. Note that, when $\sigma=1$, the
most B-robust estimator over the class of $\min\mathfrak{R}_{\alpha}%
$-estimators of $m$ is associated to $\alpha=0.4999836$, regardless the value
of $m$.

The $\psi$-function of a $\min\mathfrak{R}_{\alpha}$-estimator of the
parameter $m$, given by the formula,
\begin{equation}
\psi_{\alpha}(x,m)=\alpha(\alpha+1)^{\frac{\alpha}{2(\alpha+1)}}\sigma
^{-\frac{3\alpha+2}{\alpha+1}}(\sqrt{2\pi})^{-\frac{\alpha}{\alpha+1}%
}(x-m)\exp\left(  -\frac{\alpha}{2}\left(  \frac{x-m}{\sigma}\right)
^{2}\right)
\end{equation}
is redescending w.r.t. $x$. Then, the asymptotic breakdown point of the
corresponding estimator is 0.5 according to the results regarding redescending
M-estimators of location parameters presented in Marona et al. (2006), p.59.

The asymptotic relative efficiency is
\begin{equation}
\mathrm{ARE}(T_{\alpha},P_{m})=\frac{(2\alpha+1)^{3/2}}{(\alpha+1)^{3}}.
\end{equation}
Efficiency calculations are presented in the third row of Table 1. Small
$\alpha$ $\min\mathfrak{R}_{\alpha}$ estimation continues to retain high
efficiency. \begin{figure}[ptb]
\begin{center}
\includegraphics[width=10cm,height=10cm]{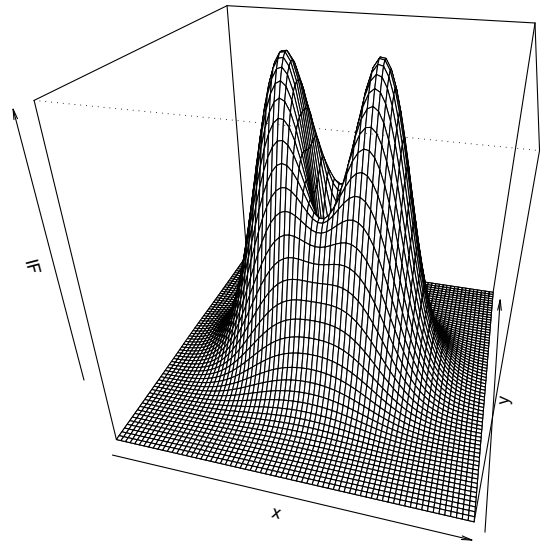}
\end{center}
\par
\vspace{-5mm}\caption{The norm of the influence function of a $\min
\mathfrak{R}_{\alpha}$-estimator of $m=(0,0)^{t}$, when $V=\mathrm{diag}%
(2,1)$}%
\label{multf}%
\end{figure}\begin{table}[ptb]
\label{Tab1} Table 1.\vspace{3mm}\newline\text{Asymptotic relative
efficiencies of the $\min\mathfrak{R}_{\alpha}$-estimators}\newline
\par%
\begin{tabular}
[c]{lllllllll}\hline
\vspace{1mm}\hspace{2.8cm}$\alpha$: & 0.00 & 0.02 & 0.05 & 0.10 & 0.20 &
0.25 & 0.5 & 1.00\\
Model &  &  &  &  &  &  &  & \\\hline
Normal $\sigma$ & 1.00000 & 0.99884 & 0.99321 & 0.97543 & 0.91922 & 0.88527 &
0.70572 & 0.43301\\
Exponential($\theta$) & 1.00000 & 0.99846 & 0.99096 & 0.96741 & 0.89412 &
0.85070 & 0.63209 & 0.33750\\
Normal $m$ & 1.00000 & 0.99942 & 0.99660 & 0.98762 & 0.95862 & 0.94060 &
0.83805 & 0.64951\\
Mean of $\mathcal{N}_{2}(m,V)$ & 1.00000 & 0.99923 & 0.99547 & 0.98353 &
0.94521 & 0.92160 & 0.79012 & 0.56250\\
Mean of $\mathcal{N}_{3}(m,V)$ & 1.00000 & 0.99903 & 0.99434 & 0.97946 &
0.93199 & 0.90297 & 0.74493 & 0.48713\\
Mean of $\mathcal{N}_{4}(m,V)$ & 1.00000 & 0.99884 & 0.99321 & 0.97541 &
0.91896 & 0.88473 & 0.70233 & 0.42187\\\hline
\end{tabular}
\end{table}

\textit{(d) Mean of multivariate normal.} The family is $\mathcal{N}_{p}%
(m,V)$. The influence function of a $\min\mathfrak{R}_{\alpha}$-estimator of
the mean vector $m$, when $V$ is known, is
\begin{equation}
\mathrm{IF}(x;T_{\alpha},P_{m})=(\sqrt{\alpha+1})^{p+2}(x-m)\exp\left(
-\frac{\alpha}{2}(x-m)^{t}V^{-1}(x-m)\right) .\label{IFm}%
\end{equation}
This is a bounded function w.r.t. $x$, meaning that all $\min\mathfrak{R}%
_{\alpha}$-estimators of the mean vector $m$ are B-robust. In Figure
\ref{multf} the norm of the influence function (\ref{IFm}), when $m=(0,0)^{t}%
$, $V=\mathrm{diag}(2,1)$ and $\alpha=0.2$, is represented.

The limiting covariance matrix of $n^{1/2}$ times the $\min\mathfrak{R}%
_{\alpha}$-estimator of $m$, when $V$ is known, can be shown to be
\begin{equation}
\left(  \frac{\alpha+1}{\sqrt{2\alpha+1}}\right)  ^{p+2}V.
\end{equation}
Then, the asymptotic relative efficiency of a $\min\mathfrak{R}_{\alpha}%
$-estimator of $m$ is
\begin{equation}
\left(  \frac{\sqrt{2\alpha+1}}{\alpha+1}\right)  ^{p+2}.
\end{equation}
Thus, one loses efficiency for increasing $p$ if $\alpha$ is kept
fixed. In Table 1 efficiencies for some values of $\alpha$ and $p$
are presented. Again, small values of $\alpha$ ensure high
efficiency and B-robustness of the estimations.

\section{$\min\mathfrak{R}_{\alpha}$-estimators in regression models}

Suppose we have i.i.d. $(p+1)$-dimensional random vectors $(X_{i},Y_{i})$,
$i=1,\dots,n$, satisfying the linear relation
\begin{equation}
\label{model}Y_{i}=\beta^{t}X_{i}+U_{i},
\end{equation}
where the $U_{i}$'s are i.i.d. with $\mathcal{N}(0,\sigma)$ and independent of
the $X_{i}$'s. $X_{i}$ and $\beta$ are $p$-dimensional column vectors with
coordinates $(X_{i1},\dots,X_{ip})$ and $(\beta_{1},\dots,\beta_{p})$,
respectively. Call $\mathbf{X}$ the $n\times p$ matrix with elements $X_{ij}$
and assume that the distribution of $X$ is not concentrated on any subspace,
i.e. $P(a^{t}X=0)<1$, for all $a\neq0$. This condition implies that the
probability that $\mathbf{X}$ has full rank tends to one when $n\to\infty$ and
holds for example if the distribution of $X$ has density.

Let $P_{\sigma}$ be the probability measure associated to a random variable
$\mathcal{N}(0,\sigma)$ and $P_{n}(\beta)$ be the empirical measure based on
the sample $U_{1},\dots,U_{n}$, where $U_{i}=Y_{i}-\beta^{t}X_{i}$,
$i=1,\dots,n$.

The $\mathfrak{R}_{\alpha}$ pseudodistance between $P_{\sigma}$ and the
probability measure $Q$ is
\begin{eqnarray}
\label{reps}\mathfrak{R}_{\alpha}(P_{\sigma},Q) &
=&\frac{1}{\alpha+1}\ln\left(
\int p^{\alpha}_{\sigma}(x)\mathrm{d}P_{\sigma}(x)\right) +\frac{1}%
{\alpha(\alpha+1)}\ln\left( \int
q^{\alpha}(x)\mathrm{d}Q(x)\right) -\nonumber\\ &&
-\frac{1}{\alpha}\ln\left( \int
p^{\alpha}_{\sigma}(x)\mathrm{d}Q(x)\right)
.\label{Renydiv}%
\end{eqnarray}

The estimators of the parameters $\beta$ and $\sigma$ are defined
by minimizing the following empirical version of the
pseudodistance (\ref{reps}),
\begin{eqnarray}
&& \mathcal{R}_{n}(P_{\sigma},P_{n}(\beta))=\nonumber\\
& =&\frac{1}{\alpha+1}\ln\int p_{\sigma}^{\alpha+1}d\lambda+\frac{1}%
{\alpha(\alpha+1)}\ln\int\left(  \frac{1}{n}\sum_{j=1}^{n}\delta_{x-U_{j}%
}\right)  ^{\alpha}\mathrm{d}P_{n}(\beta)-\frac{1}{\alpha}\ln\int p_{\sigma
}^{\alpha}\mathrm{d}P_{n}(\beta)\nonumber\\
& =&\frac{1}{\alpha+1}\ln\int p_{\sigma}^{\alpha+1}d\lambda+\frac{1}%
{\alpha(\alpha+1)}\ln\left[  \frac{1}{n}\sum_{i=1}^{n}\left(  \frac{1}{n}%
\sum_{j=1}^{n}\delta_{U_{i}-U_{j}}\right)  ^{\alpha}\right]
-\frac{1}{\alpha }\ln\left[
\frac{1}{n}\sum_{i=1}^{n}p_{\sigma}^{\alpha}(U_{i})\right]
\nonumber\\ & =&\frac{1}{\alpha+1}\ln\int
p_{\sigma}^{\alpha+1}d\lambda+\frac{1}{\alpha
+1}\ln\left(  \frac{1}{n}\right)  -\frac{1}{\alpha}\ln\left[  \frac{1}{n}%
\sum_{i=1}^{n}p_{\sigma}^{\alpha}(U_{i})\right],\label{empiric}%
\end{eqnarray}
obtained by replacing in (\ref{Renydiv}) $Q$ with $P_{n}(\beta)$ and $q(x)$
with $\widehat{q}(x)=\frac{1}{n}\sum_{j=1}^{n}\delta_{x-U_{j}}$.

Since the middle term in the above display does not depend on $\beta$ or
$\sigma$, the estimators $\widehat{\beta}$ and $\widehat{\sigma}$ are defined
by
\begin{equation}
\arg\inf_{\beta,\sigma}\left[ \frac{1}{\alpha+1}\ln\int p_{\sigma}^{\alpha
+1}d\lambda-\frac{1}{\alpha}\ln\left(  \frac{1}{n} \sum_{i=1}^{n}p_{\sigma
}^{\alpha}(Y_{i}-\beta^{t}X_{i})\right) \right]
\end{equation}
or equivalently as
\begin{equation}
\label{eqq1}\arg\sup_{\beta,\sigma}\sum_{i=1}^{n}\frac{p_{\sigma}^{\alpha
}(Y_{i}-\beta^{t}X_{i})}{\left[ \int p_{\sigma}^{\alpha+1}d\lambda\right]
^{\frac{\alpha}{\alpha+1}}}.
\end{equation}

A simple calculation shows that
\begin{equation}
\left[ \int p_{\sigma}^{\alpha+1}d\lambda\right] ^{\frac{\alpha}{\alpha+1}%
}=(\sigma\sqrt{2\pi})^{-\frac{\alpha^{2}}{\alpha+1}}(\sqrt{\alpha+1}%
)^{-\frac{\alpha}{\alpha+1}}%
\end{equation}
and (\ref{eqq1}) writes as
\begin{equation}
\arg\sup_{\beta,\sigma}\sum_{i=1}^{n}\sigma^{-\frac{\alpha}{\alpha+1}}%
\exp\left( -\frac{\alpha}{2}\left( \frac{Y_{i}-\beta^{t}X_{i}}{\sigma}\right)
^{2}\right) .
\end{equation}

Derivating with respect to $\beta$ and $\sigma$, we see that the estimators
$\widehat{\beta}$ and $\widehat{\sigma}$ are solutions of the system:
\begin{eqnarray}
&& \sum_{i=1}^{n}\exp\left( -\frac{\alpha}{2}\left( \frac{Y_{i}-\beta^{t}X_{i}%
}{\sigma}\right) ^{2}\right) \left( \frac{Y_{i}-\beta^{t}X_{i}}{\sigma}\right)
X_{i}=0\label{eqw1}\\
&& \sum_{i=1}^{n}\exp\left( -\frac{\alpha}{2}\left( \frac{Y_{i}-\beta^{t}X_{i}%
}{\sigma}\right) ^{2}\right) \left[ \left( \frac{Y_{i}-\beta^{t}X_{i}}{\sigma
}\right) ^{2}-\frac{1}{\alpha+1}\right] =0.\label{eqq2}%
\end{eqnarray}

Note that, for $\alpha=0$, the above system corresponds to the system that
define the least square estimators of $\beta$ and $\sigma$.

The system formed by (\ref{eqw1}) and (\ref{eqq2}) can be written as
\begin{equation}
\sum_{i=1}^{n}\Psi(Z_{i},\xi)=0,
\end{equation}
where $Z_{i}=(X_{i},Y_{i})$, $\xi$ is the ($p+1$)-dimensional vector with
coordinates $(\beta,\sigma)$ and
\begin{equation}
\label{grandpsi}\Psi(Z_{i},\xi)=\left( \phi\left( \frac{Y_{i}-\beta^{t}X_{i}%
}{\sigma}\right) X_{i},\chi\left( \frac{Y_{i}-\beta^{t}X_{i}}{\sigma}\right)
\right) ^{t}%
\end{equation}
with
\begin{equation}
\phi(u)=\exp\left( -\frac{\alpha}{2}u^{2}\right) u\;\;\text{and}%
\;\;\chi(u)=\left[ u^{2}-\frac{1}{\alpha+1}\right] \exp\left( -\frac{\alpha
}{2}u^{2}\right) .
\end{equation}

\begin{figure}[ptb]
\begin{center}
\includegraphics[width=9cm,height=9cm]{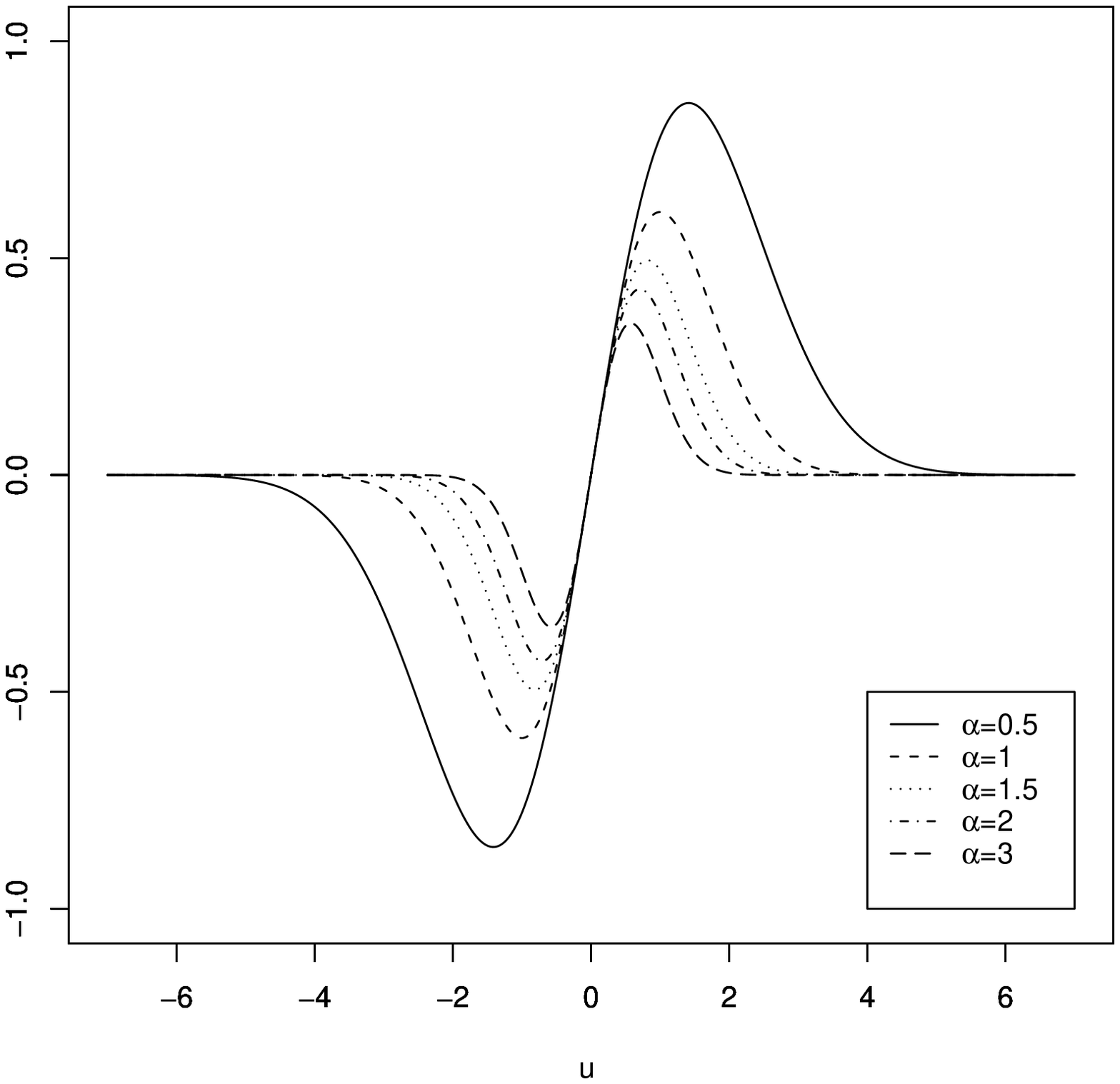}
\end{center}
\par
\vspace{-5mm}\caption{The function $\phi(u)$ corresponding to the regression
M-estimators, for different values of $\alpha$.}%
\label{redecpsi}%
\end{figure}

\begin{figure}[ptb]
\begin{center}
\includegraphics[width=9cm,height=9cm]{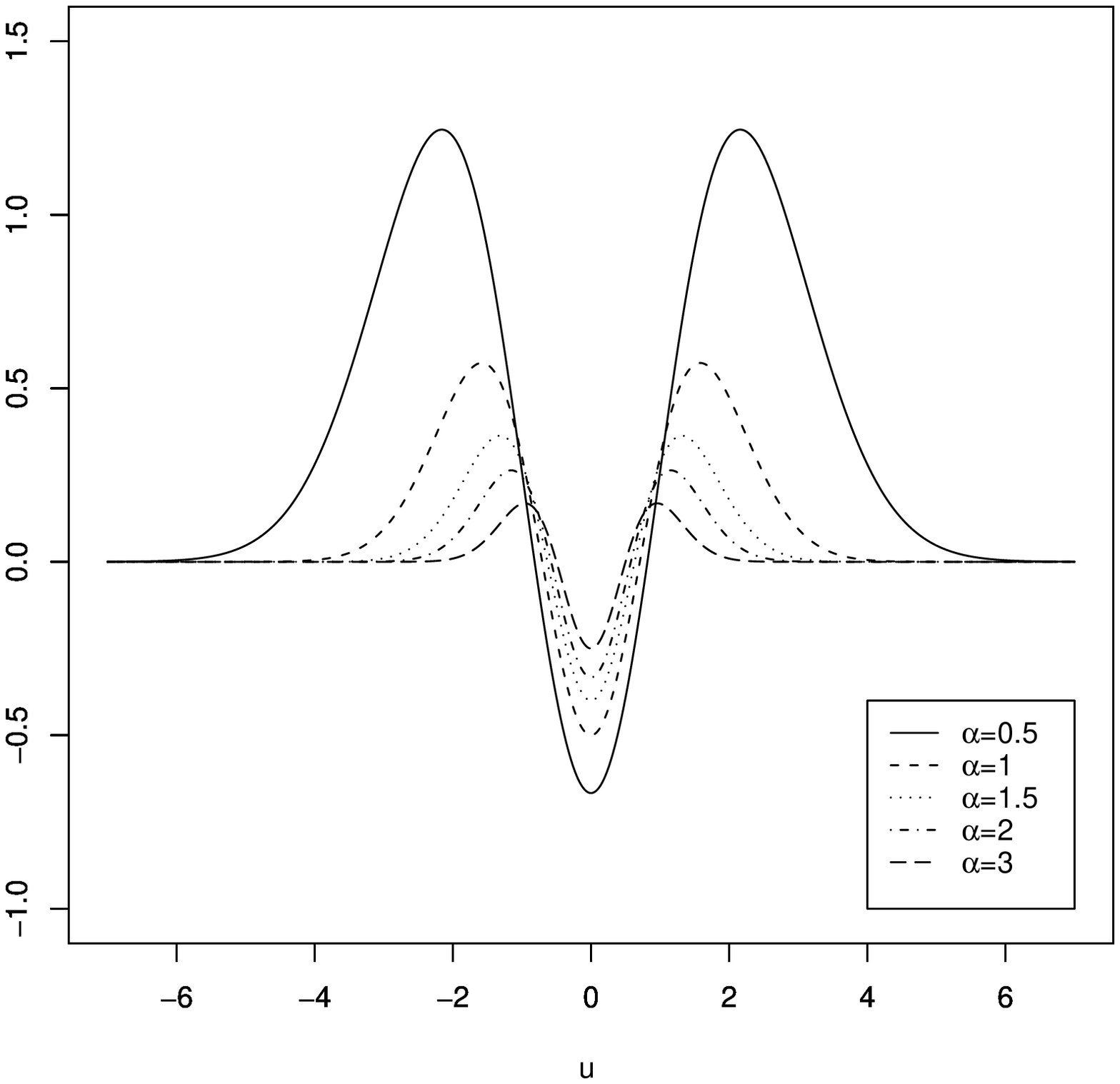}
\end{center}
\par
\vspace{-5mm}\caption{The function $\chi(u)$ corresponding to the regression
M-estimators, for different values of $\alpha$.}%
\label{chirede}%
\end{figure}

The redescending nature of the functions $\phi$ and $\chi$ can be seen in
Figure \ref{redecpsi} and in Figure \ref{chirede}.

Let $\widehat{\xi}=(\widehat{\beta},\widehat{\sigma})$. The asymptotic
normality of the M-estimator $\widehat{\xi}$ can be established by using
similar conditions with those from Theorem \ref{thasympt}. Such conditions are
satisfied by the function
\begin{equation}
h(z,\xi):=\sigma^{-\frac{\alpha}{\alpha+1}}\exp\left( -\frac{\alpha}{2}\left(
\frac{y-\beta^{t}x}{\sigma}\right) ^{2}\right)
\end{equation}
associated to the M-estimator $\widehat{\xi}$, reason for which we obtain
\begin{equation}
\sqrt{n}(\widehat{\xi}-\xi)\to\mathcal{N}_{p+1}(0,S^{-1}M(S^{-1})^{t})
\end{equation}
where
\begin{equation}
M=E\Psi(Z,\xi)\Psi(Z,\xi)^{t}\;\;\text{and}\;\;S=E\dot{\Psi}(Z,\xi),
\end{equation}
$\dot{\Psi}$ being the matrix with entries $\dot{\Psi}_{jk}:=\frac
{\partial\Psi_{j}}{\partial\xi_{k}}$.

After some calculations we find that the matrices $M$ and $S$ are
\begin{equation}
M=%
\begin{pmatrix}
\frac{\sigma^{2}}{(2\alpha+1)^{3/2}}V_{X} & 0\\
0 & \frac{\sigma^{2}(3\alpha^{2}+4\alpha+2)}{(2\alpha+1)^{5/2}(\alpha+1)^{2}}%
\end{pmatrix}
\end{equation}
and respectively
\begin{equation}
S=-%
\begin{pmatrix}
\frac{1}{(\alpha+1)^{3}}V_{X} & 0\\
0 & \frac{2}{(\alpha+1)^{5/2}}%
\end{pmatrix}
\end{equation}
where $V_{X}=EXX^{T}$.

Thus  $\widehat{\xi}$ is asymptotically normal distributed with the asymptotic
covariance matrix
\begin{equation}
\sigma^{2}%
\begin{pmatrix}
\frac{(\alpha+1)^{3}}{(2\alpha+1)^{3/2}}V_{X}^{-1} & 0\\
0 & \frac{(\alpha+1)^{3}(3\alpha^{2}+4\alpha+2)}{4(2\alpha+1)^{5/2}}%
\end{pmatrix}
.
\end{equation}

It follows that $\widehat{\beta}$ and $\widehat{\sigma}$ are asymptotically
independent and the asymptotic covariance matrix of $\widehat{\beta}$ is
$\sigma^{2}\frac{(\alpha+1)^{3}}{(2\alpha+1)^{3/2}}V_{X}^{-1}$.

Denote by $T$ and $S$ the statistical functionals corresponding to
the estimators $\widehat{\beta}$ and $\widehat{\sigma}$,
respectively. For a given probability measure $P$, these
functionals are defined through the solutions of the system
\begin{equation}
\int\Psi(z,T(P),S(P))\mathrm{d}P=0,\label{systfunc}%
\end{equation}
where
\begin{equation}
\Psi(z,\xi)=\left(  \phi\left(  \frac{y-\beta^{t}x}{\sigma}\right)
x,\chi\left(  \frac{y-\beta^{t}x}{\sigma}\right)  \right).
\end{equation}

The influence functions of the functionals $T$ and $S$ are given in the
following theorem:

\begin{theorem}
The influence functions of the functionals associated to the $\min
\mathfrak{R}_{\alpha}$-estimators of $\beta$ and $\sigma$ are
\begin{eqnarray*}
\mathrm{IF}(x_{0},y_{0};T,P_{\xi}) &
=&\sigma(\alpha+1)^{3/2}\exp\left( -\frac{\alpha}{2}\left(
\frac{y_{0}-\beta^{t}x_{0}}{\sigma}\right) ^{2}\right) \left(
\frac{y_{0}-\beta^{t}x_{0}}{\sigma}\right) V_{X}^{-1}x_{0}\\
\mathrm{IF}(x_{0},y_{0};S,P_{\xi}) &
=&\frac{(\alpha+1)^{5/2}}{2}\exp\left( -\frac{\alpha}{2}\left(
\frac{y_{0}-\beta^{t}x_{0}}{\sigma}\right) ^{2}\right)
\left[ \left( \frac{y_{0}-\beta^{t}x_{0}}{\sigma}\right) ^{2}-\frac{1}%
{\alpha+1}\right],
\end{eqnarray*}
$P_{\xi}$ being the probability measure associated to $Z$.
\end{theorem}

\textit{Proof.} The system (\ref{systfunc}) can be written as
\begin{eqnarray}
& \int\phi\left( \frac{y-T(P)^{t}x}{S(P)}\right) x\mathrm{d}%
P(x,y)=0\nonumber\\
& \int\chi\left( \frac{y-T(P)^{t}x}{S(P)}\right) \mathrm{d}%
P(x,y)=0.\label{systfunc2}
\end{eqnarray}

We consider the contaminated model $\widetilde{P}_{\varepsilon, x_{0},y_{0}%
}=(1-\varepsilon)P_{\xi}+\varepsilon\delta_{(x_{0},y_{0})}$, where
$(x_{0},y_{0})$ is an arbitrary point from $\mathbb{R}^{p}\times\mathbb{R}$.
For this model, the system (\ref{systfunc2}) writes as
\begin{eqnarray*}
&& (1-\varepsilon)\int\phi\left(
\frac{y-T(\widetilde{P}_{\varepsilon,
x_{0},y_{0}})^{t}x}{S(\widetilde{P}_{\varepsilon,
x_{0},y_{0}})}\right)
x\mathrm{d}P_{\xi}(x,y)+\varepsilon\phi\left(
\frac{y_{0}-T(\widetilde {P}_{\varepsilon,
x_{0},y_{0}})^{t}x_{0}}{S(\widetilde{P}_{\varepsilon,
x_{0},y_{0}})}\right) x_{0}=0\\ && (1-\varepsilon)\int\chi\left(
\frac{y-T(\widetilde{P}_{\varepsilon,
x_{0},y_{0}})^{t}x}{S(\widetilde{P}_{\varepsilon,
x_{0},y_{0}})}\right) \mathrm{d}P_{\xi}(x,y)+\varepsilon\chi\left(
\frac{y_{0}-T(\widetilde {P}_{\varepsilon,
x_{0},y_{0}})^{t}x_{0}}{S(\widetilde{P}_{\varepsilon,
x_{0},y_{0}})}\right) =0.
\end{eqnarray*}
Derivating with respect to $\varepsilon$ and taking the derivative in
$\varepsilon=0$, after some calculations, we find
\begin{eqnarray*}
\mathrm{IF}(x_{0},y_{0};T,P_{\xi}) &
=&\sigma(\alpha+1)^{3/2}\phi\left(
\frac{y_{0}-\beta^{t}x_{0}}{\sigma}\right) V_{X}^{-1}x_{0}\\
& =&\sigma(\alpha+1)^{3/2}\exp\left( -\frac{\alpha}{2}\left( \frac{y_{0}%
-\beta^{t}x_{0}}{\sigma}\right) ^{2}\right) \left( \frac{y_{0}-\beta^{t}x_{0}%
}{\sigma}\right) V_{X}^{-1}x_{0}%
\end{eqnarray*}
and
\begin{eqnarray*}
\mathrm{IF}(x_{0},y_{0};S,P_{\xi}) &
=&\frac{(\alpha+1)^{5/2}}{2}\chi\left(
\frac{y_{0}-\beta^{t}x_{0}}{\sigma}\right) \\ & =&
\frac{(\alpha+1)^{5/2}}{2}\exp\left( -\frac{\alpha}{2}\left( \frac
{y_{0}-\beta^{t}x_{0}}{\sigma}\right) ^{2}\right) \left[ \left(
\frac {y_{0}-\beta^{t}x_{0}}{\sigma}\right)
^{2}-\frac{1}{\alpha+1}\right].
\end{eqnarray*}

Since $\chi$ is redescending, the estimator $\widehat{\sigma}$ has the
influence function bounded and hence is B-robust. On the other hand,
$\mathrm{IF}(x_{0},y_{0};T,P_{\xi})$ will tend to infinity only when $x_{0}$
tends to infinity and $\left| \frac{y_{0}-\beta^{t}x_{0}}{\sigma}\right| \leq
k$, for some $k$. This means that large outliers have no influence on the estimates.

\section{Simulation results}

In this section we present some simulation studies in order to illustrate the
performance of $\min\mathfrak{R}_{\alpha}$-estimators in finite samples.

First, we considered the scale normal model with known mean. We estimated the
scale parameter $\sigma$ by using the $\min\mathfrak{R}_{\alpha}$-estimator
which is obtained as solution of the equation
\begin{equation}
\sum_{i=1}^{n}\left[ \left( \frac{X_{i}-m}{\sigma}\right) ^{2}-\frac{1}%
{\alpha+1}\right] \exp\left( -\frac{\alpha}{2}\left( \frac{X_{i}-m}{\sigma
}\right) ^{2}\right)  =0,
\end{equation}
$m$ being the known mean.

To make some comparisons, we also considered the minimum density power
divergence estimator of Basu et al. (1998) (in the present paper we will
denote it by $\min\mathfrak{D}_{\alpha}$-estimator). For the scale normal
model, this estimator is solution of the equation
\begin{equation}
\frac{\alpha}{(\alpha+1)^{3/2}}+\frac{1}{n}\sum_{i=1}^{n}\left[ \left(
\frac{X_{i}-m}{\sigma}\right) ^{2}-1\right] \exp\left( -\frac{\alpha}{2}\left(
\frac{X_{i}-m}{\sigma}\right) ^{2}\right)  =0,
\end{equation}

In a first Monte Carlo experiment, 5000 samples of size $n=100$ were generated
from the scale normal model $\mathcal{N}(0,1)$ with mean $m=0$ known,
$\sigma=1$ being the parameter to be estimated. In a second Monte Carlo
experiment we generated 5000 samples with $100$ observations, for each sample
95 observations being generated from $\mathcal{N}(0,1)$ and 5 from
$\mathcal{N}(2,1)$, and then we generated 5000 samples with 100 observations,
for each sample 90 observations being generated from $\mathcal{N}(0,1)$ and 10
from $\mathcal{N}(2,1)$. For each sample we computed $\min\mathfrak{R}%
_{\alpha}$-estimators and $\min\mathfrak{D}_{\alpha}$-estimators corresponding
to $\alpha\in\{0.02,0.05,0.1,0.2,0.25,0.5,1\}$ and the MLE for $\alpha=0$.

In Table 2 we present the mean estimated scale $\widehat{\sigma}$ and
simulation based estimates of the MSE defined by
\begin{equation}
\widehat{\mathrm{MSE}}:=\frac{1}{n_{s}}\sum_{i=1}^{n_{s}}(\widehat{\sigma}%
_{i}-\sigma)^{2}%
\end{equation}
where $n_{s}$ denotes the number of samples (5000 in our study) and
$\widehat{\sigma}_{i}$ represents an estimate of $\sigma=1$ obtained on the
basis of the $i$th sample.

As it can be seen, both the $\min\mathfrak{R}_{\alpha}$-estimators and
$\min\mathfrak{D}_{\alpha}$-estimators perform well under the model. Under
contamination, the $\min\mathfrak{R}_{\alpha}$-estimator with $\alpha=1$ gives
the best results in terms of robustness, while keeping small empirical MSE.
However, the $\min\mathfrak{R}_{\alpha}$-estimators, as well as the
$\min\mathfrak{D}_{\alpha}$-estimators, exhibit outlier resistance properties
even for small values of $\alpha$. For example, in the case of $5\%$
contamination, the estimates of $\sigma=1$ obtained for $\alpha=0.2$ are
1.07494, respectively 1.07562, fairly close to the estimates obtained for
$\alpha=1$. In this case, the $\min\mathfrak{R}_{\alpha}$-estimator combines
robustness with the asymptotic relative efficiency 0.91922.

Similar results are presented in Table 3 and Table 4, where 5$\%$ or 10$\%$
from data come from the contaminating distribution $\mathcal{N}(0,3)$ or from
$\delta_{10}$. When the contamining distribution is $\delta_{10}$, the
$\min\mathfrak{R}_{\alpha}$-estimators have strong robustness properties, and
this can also be explained by the influence function which is redescending, as
it can be seen in Figure \ref{Gr1}.

In the second example, our estimation method is applied to the location normal
model $\mathcal{N}(0,1)$, $\sigma=1$ being known. Here we compute the
$\min\mathfrak{R}_{\alpha}$-estimates as solutions of equation
\begin{equation}
\sum_{i=1}^{n}\left( X_{i}-m\right) \exp\left( -\frac{\alpha}{2}\left(
\frac{X_{i}-m}{\sigma}\right) ^{2}\right)  =0.
\end{equation}
We consider the case of no outliers and the cases of $5\%$ or $10\%$ outliers
coming from the model $\mathcal{N}(2,1)$. The results are given in Table 5.
Again, the choice $\alpha=0.2$ provides robustness and high efficiency of the
estimation procedure. When the outliers come from $\delta_{10}$, we obtain
very good results in terms of robustness, even for very small values of
$\alpha$, as it can be seen in Table 6. These results are in accordance with
the redescending nature of the influence functions represented in Figure
\ref{Gr3}.

Our examples show that increasing $\alpha$ leads to estimators
which are far more robust than the maximum likelihood estimator.
The simulation results suggest that $\alpha$ between 0.1 and 0.25
provides competitive estimators in terms of robustness and
efficiency.

\begin{center}
\begin{table}[ptb]
\label{table2} Table 2.\vspace{3mm}\newline\text{Simulation results for
$\min\mathfrak{R}_{\alpha}$-estimators, $\min\mathfrak{D}_{\alpha}$-estimators
and MLE of the }\newline\text{parameter $\sigma=1$ when data are generated
from the model $\mathcal{N}(0,1)$, when 95 data}\newline\text{are generated
from the model $\mathcal{N}(0,1)$ and 5 data from $\mathcal{N}(2,1)$,
respectively when }\newline\vspace{2mm}\text{90 data are generated from the
model $\mathcal{N}(0,1)$ and 10 data from $\mathcal{N}(2,1)$.}\newline%
\begin{tabular}
[c]{lllllllll}\hline
& \!\!no outliers\! &  &  & \!\!5\% outliers\! &  &  & \!\!10\% outliers\! &
\\\cline{2-3}\cline{5-6}\cline{8-9}%
\vspace{1mm} & $\widehat{\sigma}$ & $\widehat{\mathrm{MSE}}$ &  &
$\widehat{\sigma}$ & $\widehat{\mathrm{MSE}}$ &  & $\widehat{\sigma}$ &
$\widehat{\mathrm{MSE}}$\\\hline
&  &  &  &  &  &  &  & \\
MLE &  &  &  &  &  &  &  & \\
$\alpha=0$ & 0.99763 & 0.00503 &  & 1.09289 & 0.01446 &  & 1.17999 & 0.03888\\
&  &  &  &  &  &  &  & \\
$\min\mathfrak{R}_{\alpha}$ &  &  &  &  &  &  &  & \\
$\alpha=0.02$ & 0.99987 & 0.00501 &  & 1.09216 & 0.01420 &  & 1.17886 &
0.03827\\
$\alpha=0.05$ & 1.00022 & 0.00504 &  & 1.08902 & 0.01357 &  & 1.17445 &
0.03663\\
$\alpha=0.1$ & 1.00069 & 0.00514 &  & 1.08398 & 0.01272 &  & 1.16712 &
0.03412\\
$\alpha=0.2$ & 1.00122 & 0.00545 &  & 1.07494 & 0.01162 &  & 1.15310 &
0.02999\\
$\alpha=0.25$ & 1.00137 & 0.00566 &  & 1.07100 & 0.01131 &  & 1.14659 &
0.02835\\
$\alpha=0.5$ & 1.00142 & 0.00710 &  & 1.05610 & 0.01122 &  & 1.11981 &
0.02348\\
$\alpha=1$ & 0.99956 & 0.01173 &  & 1.03931 & 0.01494 &  & 1.08746 & 0.02315\\
&  &  &  &  &  &  &  & \\
$\min\mathfrak{D}_{\alpha}$ &  &  &  &  &  &  &  & \\
$\alpha=0.02$ & 0.99977 & 0.00463 &  & 1.09233 & 0.01398 &  & 1.17940 &
0.03828\\
$\alpha=0.05$ & 1.00012 & 0.00467 &  & 1.08926 & 0.01336 &  & 1.17505 &
0.03665\\
$\alpha=0.1$ & 1.00060 & 0.00477 &  & 1.08434 & 0.01252 &  & 1.16784 &
0.03415\\
$\alpha=0.2$ & 1.00123 & 0.00508 &  & 1.07562 & 0.01144 &  & 1.15420 &
0.03007\\
$\alpha=0.25$ & 1.00146 & 0.00527 &  & 1.07194 & 0.01113 &  & 1.14802 &
0.02848\\
$\alpha=0.5$ & 1.00217 & 0.00648 &  & 1.05945 & 0.01093 &  & 1.12494 &
0.02390\\
$\alpha=1$ & 1.00326 & 0.00882 &  & 1.05181 & 0.01251 &  & 1.10781 &
0.02259\\\hline
\end{tabular}
\end{table}

\begin{table}[ptb]
Table 3.\vspace{3mm}\newline\text{Simulation results for $\min\mathfrak{R}%
_{\alpha}$-estimators, $\min\mathfrak{D }_{\alpha}$-estimators and MLE of the
}\newline\text{parameter $\sigma=1$ when data are generated from the model
$\mathcal{N}(0,1)$, when 95 data}\newline\text{are generated from the model
$\mathcal{N}(0,1)$ and 5 data from $\mathcal{N}(0,3)$, respectively when
}\newline\vspace{2mm}\text{90 data are generated from the model $\mathcal{N}%
(0,1)$ and 10 data from $\mathcal{N}(0,3)$.}\newline%
\begin{tabular}
[c]{lllllllll}\hline
& \!\!no outliers\! &  &  & \!\!5\% outliers\! &  &  & \!\!10\% outliers\! &
\\\cline{2-3}\cline{5-6}\cline{8-9}%
\vspace{1mm} & $\widehat{\sigma}$ & $\widehat{\mathrm{MSE}}$ &  &
$\widehat{\sigma}$ & $\widehat{\mathrm{MSE}}$ &  & $\widehat{\sigma}$ &
$\widehat{\mathrm{MSE}}$\\\hline
&  &  &  &  &  &  &  & \\
MLE &  &  &  &  &  &  &  & \\
$\alpha=0$ & 0.99794 & 0.00498 &  & 1.17726 & 0.04887 &  & 1.33251 & 0.13507\\
&  &  &  &  &  &  &  & \\
$\min\mathfrak{R}_{\alpha}$ &  &  &  &  &  &  &  & \\
$\alpha=0.02$ & 0.99749 & 0.00493 &  & 1.15713 & 0.03788 &  & 1.30542 &
0.11284\\
$\alpha=0.05$ & 0.99783 & 0.00496 &  & 1.13024 & 0.02669 &  & 1.26450 &
0.08485\\
$\alpha=0.1$ & 0.99827 & 0.00505 &  & 1.09683 & 0.01692 &  & 1.20663 &
0.05362\\
$\alpha=0.2$ & 0.99874 & 0.00536 &  & 1.06176 & 0.01059 &  & 1.13522 &
0.02695\\
$\alpha=0.25$ & 0.99884 & 0.00557 &  & 1.05226 & 0.00955 &  & 1.11441 &
0.02152\\
$\alpha=0.5$ & 0.99869 & 0.00709 &  & 1.03035 & 0.00905 &  & 1.06610 &
0.01367\\
$\alpha=1$ & 0.99670 & 0.01198 &  & 1.01659 & 0.01356 &  & 1.03857 & 0.01598\\
&  &  &  &  &  &  &  & \\
$\min\mathfrak{D}_{\alpha}$ &  &  &  &  &  &  &  & \\
$\alpha=0.02$ & 0.99870 & 0.00487 &  & 1.15614 & 0.03690 &  & 1.30385 &
0.11189\\
$\alpha=0.05$ & 0.99914 & 0.00489 &  & 1.12922 & 0.02596 &  & 1.26309 &
0.08411\\
$\alpha=0.1$ & 0.99973 & 0.00497 &  & 1.09665 & 0.01677 &  & 1.20657 &
0.05379\\
$\alpha=0.2$ & 1.00051 & 0.00525 &  & 1.06368 & 0.01083 &  & 1.13856 &
0.02810\\
$\alpha=0.25$ & 1.00078 & 0.00543 &  & 1.05509 & 0.00982 &  & 1.11932 &
0.02285\\
$\alpha=0.5$ & 1.00143 & 0.00659 &  & 1.03786 & 0.00906 &  & 1.07956 &
0.01526\\
$\alpha=1$ & 1.00203 & 0.00891 &  & 1.03484 & 0.01101 &  & 1.07089 &
0.01588\\\hline
\end{tabular}
\end{table}

\begin{table}[ptb]
Table 4.\vspace{3mm}\newline\text{Simulation results for $\min\mathfrak{R}%
_{\alpha}$-estimators, $\min\mathfrak{D}_{\alpha}$-estimators}\newline%
\text{and MLE of the parameter $\sigma=1$ when the data are generated}%
\newline\text{from the model $\mathcal{N}(0,1)$, when 95 data are generated
from the }\newline\text{model $\mathcal{N}(0,1)$ and 5 data from $\delta_{10}%
$.}\newline
\par%
\begin{tabular}
[c]{llllll}\hline
& \!\!no outliers\! &  &  & \!\!5\% outliers\! & \\\cline{2-3}\cline{5-6}%
\vspace{1mm} & $\widehat{\sigma}$ & $\widehat{\mathrm{MSE}}$ &  &
$\widehat{\sigma}$ & $\widehat{\mathrm{MSE}}$\\\hline
&  &  &  &  & \\
MLE &  &  &  &  & \\
$\alpha=0$ & 0.99859 & 0.00497 &  & 2.43937 & 2.07260\\
&  &  &  &  & \\
$\min\mathfrak{R}_{\alpha}$ &  &  &  &  & \\
$\alpha=0.02$ & 0.99731 & 0.00494 &  & 2.28640 & 1.65613\\
$\alpha=0.05$ & 0.99773 & 0.00497 &  & 1.95068 & 0.90877\\
$\alpha=0.1$ & 0.99828 & 0.00506 &  & 1.03369 & 0.01115\\
$\alpha=0.2$ & 0.99900 & 0.00538 &  & 0.99922 & 0.00575\\
$\alpha=0.25$ & 0.99923 & 0.00559 &  & 0.99913 & 0.00592\\
$\alpha=0.5$ & 0.99968 & 0.00708 &  & 0.99947 & 0.00747\\
$\alpha=1$ & 0.99875 & 0.01176 &  & 0.99832 & 0.01238\\
&  &  &  &  & \\
$\min\mathfrak{D}_{\alpha}$ &  &  &  &  & \\
$\alpha=0.02$ & 0.99567 & 0.00506 &  & 2.28566 & 1.65427\\
$\alpha=0.05$ & 0.99602 & 0.00509 &  & 1.94971 & 0.90702\\
$\alpha=0.1$ & 0.99646 & 0.00517 &  & 1.03533 & 0.01155\\
$\alpha=0.2$ & 0.99700 & 0.00544 &  & 1.00358 & 0.00585\\
$\alpha=0.25$ & 0.99716 & 0.00563 &  & 1.00522 & 0.00602\\
$\alpha=0.5$ & 0.99755 & 0.00678 &  & 1.01521 & 0.00748\\
$\alpha=1$ & 0.99818 & 0.00905 &  & 1.03344 & 0.01082\\\hline
\end{tabular}
\end{table}

\begin{table}[ptb]
Table 5.\vspace{3mm}\newline\text{Simulation results for $\min\mathfrak{R}%
_{\alpha}$-estimators and MLE of the parameter $m=0$ when}\newline\text{data
are generated from the model $\mathcal{N}(0,1)$, when 95 data are generated
from the }\newline\text{model $\mathcal{N}(0,1)$ and 5 data from
$\mathcal{N}(2,1)$, respectively when 90 data are generated }\newline%
\vspace{2mm}\text{from the model $\mathcal{N}(0,1)$ and 10 data from
$\mathcal{N}(2,1)$.}\newline%
\begin{tabular}
[c]{lllllllll}\hline
& \!\!no outliers\! &  &  & \!\!5\% outliers\! &  &  & \!\!10\% outliers\! &
\\\cline{2-3}\cline{5-6}\cline{8-9}%
\vspace{1mm} & $\widehat{m}$ & $\widehat{\mathrm{MSE}}$ &  & $\widehat{m}$ &
$\widehat{\mathrm{MSE}}$ &  & $\widehat{m}$ & $\widehat{\mathrm{MSE}}$\\\hline
&  &  &  &  &  &  &  & \\
MLE &  &  &  &  &  &  &  & \\
$\alpha=0$ & 0.00158 & 0.01004 &  & 0.10161 & 0.02054 &  & 0.20116 & 0.05063\\
&  &  &  &  &  &  &  & \\
$\min\mathfrak{R}_{\alpha}$ &  &  &  &  &  &  &  & \\
$\alpha=0.02$ & 0.00075 & 0.01003 &  & 0.09723 & 0.01956 &  & 0.19573 &
0.04853\\
$\alpha=0.05$ & 0.00072 & 0.01006 &  & 0.09257 & 0.01873 &  & 0.18769 &
0.04554\\
$\alpha=0.1$ & 0.00067 & 0.01015 &  & 0.08568 & 0.01768 &  & 0.17557 &
0.04141\\
$\alpha=0.2$ & 0.00059 & 0.01045 &  & 0.07457 & 0.01640 &  & 0.15539 &
0.03549\\
$\alpha=0.25$ & 0.00055 & 0.01065 &  & 0.07007 & 0.01605 &  & 0.14698 &
0.03339\\
$\alpha=0.5$ & 0.00036 & 0.01194 &  & 0.05427 & 0.01583 &  & 0.11651 &
0.02774\\
$\alpha=1$ & 0.00004 & 0.01550 &  & 0.03900 & 0.01864 &  & 0.08580 &
0.02675\\\hline
\end{tabular}
\end{table}

\begin{table}[ptb]
Table 6.\vspace{3mm}\newline\text{Simulation results for $\min\mathfrak{R}%
_{\alpha}$-estimators and MLE of the parameter $m=0$ when}\newline\text{data
are generated from the model $\mathcal{N}(0,1)$, when 95 data are generated
from the }\newline\text{model $\mathcal{N}(0,1)$ and 5 data from $\delta_{10}%
$, respectively when 90 data are generated from}\newline\vspace{2mm}\text{the
model $\mathcal{N}(0,1)$ and 10 data from $\delta_{10}$.}\newline%
\begin{tabular}
[c]{lllllllll}\hline
& \!\!no outliers\! &  &  & \!\!5\% outliers\! &  &  & \!\!10\% outliers\! &
\\\cline{2-3}\cline{5-6}\cline{8-9}%
\vspace{1mm} & $\widehat{m}$ & $\widehat{\mathrm{MSE}}$ &  & $\widehat{m}$ &
$\widehat{\mathrm{MSE}}$ &  & $\widehat{m}$ & $\widehat{\mathrm{MSE}}$\\\hline
&  &  &  &  &  &  &  & \\
MLE &  &  &  &  &  &  &  & \\
$\alpha=0$ & 0.00227 & 0.00999 &  & 0.50186 & 0.26125 &  & 1.00183 & 1.01247\\
&  &  &  &  &  &  &  & \\
$\min\mathfrak{R}_{\alpha}$ &  &  &  &  &  &  &  & \\
$\alpha=0.02$ & 0.00018 & 0.01021 &  & 0.20421 & 0.05278 &  & 0.44055 &
0.20618\\
$\alpha=0.05$ & 0.00023 & 0.01023 &  & 0.04819 & 0.01342 &  & 0.10345 &
0.02288\\
$\alpha=0.1$ & 0.00033 & 0.01033 &  & 0.00499 & 0.01088 &  & 0.00999 &
0.01159\\
$\alpha=0.2$ & 0.00054 & 0.01064 &  & 0.00106 & 0.01111 &  & 0.00147 &
0.01166\\
$\alpha=0.25$ & 0.00065 & 0.01085 &  & 0.00112 & 0.01132 &  & 0.00150 &
0.01188\\
$\alpha=0.5$ & 0.00118 & 0.01219 &  & 0.00158 & 0.01275 &  & 0.00194 &
0.01336\\
$\alpha=1$ & 0.00198 & 0.01582 &  & 0.00230 & 0.01658 &  & 0.00271 &
0.01740\\\hline
\end{tabular}
\end{table}
\end{center}

\end{document}